\documentclass{article}
\usepackage[scale=0.8]{geometry}

\usepackage{amsfonts}
\usepackage{amsthm}
\usepackage{graphicx}
\usepackage{mathtools}
\usepackage{algorithmic}
\usepackage{bm}
\usepackage{subcaption}
\usepackage{booktabs}
\usepackage{tikz}
\usetikzlibrary{arrows.meta, calc, quotes, tikzmark, shapes.misc, positioning, arrows, decorations.pathreplacing}
\tikzset{cross/.style={cross out, draw=black, minimum size=2*(#1-\pgflinewidth), inner sep=0pt, outer sep=0pt},
    cross/.default={1pt}}
\usepackage{pgfplots}
\pgfplotsset{compat=1.18}
\usepackage{xspace}
\usepackage{enumitem}
\usepackage{authblk}

\usepackage[colorlinks=true]{hyperref}
\usepackage[capitalise,nameinlink,noabbrev]{cleveref}

\newtheorem{theorem}{Theorem}[section]
\newtheorem{lemma}[theorem]{Lemma}

\theoremstyle{definition}
\newtheorem{definition}[theorem]{Definition}

\theoremstyle{remark}
\newtheorem{remark}[theorem]{Remark}

\numberwithin{equation}{section}

\DeclareMathOperator*{\equals}{=}

\crefname{corollary}{Corollary}{Corollaries}

\makeatletter
\newsavebox\myboxA
\newsavebox\myboxB
\newlength\mylenA

\newcommand*\xoverline[2][0.75]{%
    \sbox{\myboxA}{$\m@th#2$}%
    \setbox\myboxB\null
    \ht\myboxB=\ht\myboxA%
    \dp\myboxB=\dp\myboxA%
    \wd\myboxB=#1\wd\myboxA
    \sbox\myboxB{$\m@th\overline{\copy\myboxB}$}
    \setlength\mylenA{\the\wd\myboxA}
    \addtolength\mylenA{-\the\wd\myboxB}%
    \ifdim\wd\myboxB<\wd\myboxA%
       \rlap{\hskip 0.5\mylenA\usebox\myboxB}{\usebox\myboxA}%
    \else
        \hskip -0.5\mylenA\rlap{\usebox\myboxA}{\hskip 0.5\mylenA\usebox\myboxB}%
    \fi}

\makeatother

\newcommand{\pd}{\mathcal{\partial}}
\newcommand{\WB}{\text{WB}}
\newcommand{\imh}{{i - \frac 1 2}}
\newcommand{\iph}{{i + \frac 1 2}}
\newcommand{\jmh}{{j - \frac 1 2}}
\newcommand{\jph}{{j + \frac 1 2}}
\newcommand{\imhinline}{{i - 1/2}}
\newcommand{\iphinline}{{i + 1/2}}
\newcommand{\dx}{{\Delta x}}
\newcommand{\dy}{{\Delta y}}
\newcommand{\dt}{{\Delta t}}
\newcommand{\normalx}{{\boldsymbol e_x}}
\newcommand{\normaly}{{\boldsymbol e_y}}

\newcommand{\ARS}{A-RS\xspace}
\newcommand{\ARSs}{A-RSs\xspace}
\newcommand{\HLL}{\text{HLL}}

\DeclareFontFamily{U}{mathx}{}
\DeclareFontShape{U}{mathx}{m}{n}{<-> mathx10}{}
\DeclareSymbolFont{mathx}{U}{mathx}{m}{n}
\DeclareMathAccent{\widehat}{0}{mathx}{"70}
\DeclareMathAccent{\widecheck}{0}{mathx}{"71}

\makeatletter
\def\widebreve{\mathpalette\wide@breve}
\def\wide@breve#1#2{\sbox\z@{$#1#2$}%
     \mathop{\vbox{\m@th\ialign{##\crcr
\kern0.08em\brevefill#1{0.8\wd\z@}\crcr\noalign{\nointerlineskip}%
                    $\hss#1#2\hss$\crcr}}}\limits\hspace*{-0.175em}}
\def\brevefill#1#2{$\m@th\sbox\tw@{$#1($}%
  \hss\resizebox{#2}{\wd\tw@}{\rotatebox[origin=c]{90}{\upshape(}}\hss$}
\makeatletter

\title{Towards a fully well-balanced and entropy-stable scheme for the Euler equations with gravity: preserving isentropic steady solutions}
\author[1]{Christophe Berthon}
\author[2]{Victor Michel-Dansac}
\author[2]{Andrea Thomann}
\affil[1]{\small Laboratoire de Mathématiques Jean Leray, CNRS UMR 6629, Université de Nantes, 2 rue de la Houssinière, BP 92208, 44322 Nantes Cedex 3, France}
\affil[2]{\small Université de Strasbourg, CNRS, Inria, IRMA, F-67000 Strasbourg, France}

\begin{document}







\date{\today}



\maketitle

\begin{abstract}
    The present work concerns the derivation of a numerical
    scheme to approximate weak solutions of the Euler equations with a
    gravitational source term. The designed scheme is proved to be
    fully well-balanced since it is able to exactly preserve all
    moving equilibrium solutions,
    as well as the corresponding steady solutions at rest obtained
    when the velocity vanishes.
    Moreover, the proposed scheme
    is entropy-preserving since it satisfies all fully discrete
    entropy inequalities. In addition, in order to satisfy the required
    admissibility of the approximate solutions, the positivity of both
    approximate density and pressure is established. Several numerical
    experiments attest the relevance of the developed numerical method.
    An extension to two-dimensional problems is given, applying the one-dimensional framework direction by direction on Cartesian grids.
\end{abstract}

\section{Introduction}
\label{sec:Intro}
The present paper is devoted to designing a numerical scheme to
approximate the weak solutions of the well-known Euler equations with a given gravitational source term \cite{Xu2010,Luo2011,XinShu2012,Vides2014,Desveaux2016,KliPupSem2019,Thomann2019}.
From a numerical point of view, such a system poses three main challenges:
\begin{align}
     & \text{the physical admissibility of the approximate solutions,} \tag{a}\label{i} \\
     & \text{satisfying the discrete entropy stability,}\tag{b}\label{ii}               \\
     & \text{the suitable capture of the steady solutions.}\tag{c} \label{iii}
\end{align}
The admissibility of the discrete solutions, namely preserving the positivity of both density and pressure, is the basic property to be satisfied by any first-order numerical scheme.
For instance, the Godunov scheme, the HLL and HLLC schemes, as well as relaxation
schemes are usual numerical methods which satisfy this essential robustness
property, see
\cite{HarLaxLee1983,Toro1994,Jin1996,Coquel1998,Chalons2008e,Toro2009},
but the list is clearly not exhaustive.
However, it is worth noting that the set of physically admissible states is not convex.
Hence, specific techniques, like entropy inequalities,
are required to get a positive approximate pressure, see for instance \cite{Berthon2006}.

Moreover, entropy inequalities are essential to rule out some
non-physical solutions. Indeed, since the system under consideration
is hyperbolic, discontinuous solutions may form in finite time \cite{Godlewski1991,Dafermos2010}. The physical admissibility
of the discontinuities are thus governed by the entropy
inequalities. The homogeneous Euler system contains a large range of
entropy inequalities, since they are defined up to a smooth function,
according to prescribed restrictions \cite{Tadmor1986}.

In this work, we are concerned with a non-homogeneous Euler system, supplemented by a
gravitational source term. It is thus very interesting to remark that
the gravitational source term is compatible with the entropies, as shown in
\cite{Desveaux2016}, since it contains the same set of entropy
functions as the usual homogeneous Euler system. From a numerical
point of view, this essential entropy stability must be preserved. The
establishment of discrete entropy inequalities is known to be very
challenging for usual homogeneous hyperbolic systems. For instance, the
Godunov scheme, the HLL and HLLC schemes, or the Suliciu and Xin-Jin relaxation schemes are known to be entropy-satisfying for the Euler equations, see~\cite{HarLaxLee1983,Godlewski1991,Bouchut2004,Berthon2006,Chalons2008e,Desveaux2016}. Here, however, we have to take into account additional source terms.
In order to recover this crucial stability property, we adopt some recent Godunov-type
techniques \cite{MicBerClaFou2016,MicBerClaFou2017} which are well-designed
to approximate source terms.
In particular, in one of our main results, we show that the specific entropy
given by the celebrated approximate Riemann solver by
Harten, Lax and van Leer~\cite{HarLaxLee1983}
satisfies a decreasing principle when a convex function composition is applied.
This valuable and useful principle is the main argument
to prove that the approximate solutions generated by
our new approximate Riemann solver satisfy the discrete entropy inequalities.

The last challenge that we need to address is posed by accurately capturing steady solutions.
Indeed, due to the source term, the system admits non-obvious steady solutions
\cite{Bouchut2004,Pares2004,Desveaux2016,FraMen2016}. These particular equilibrium states
are of strong interest from a numerical point of view. Indeed, huge
errors can be introduced by the source term discretization, and resulting spurious waves
can be propagated over the whole computational domain.
This numerical failure has been underlined a long time ago in
\cite{Bermudez1994,Cargo1994,Greenberg1996,Gosse2000} for
shallow water simulations. For over two decades, much work was
devoted to the derivation of schemes preserving steady solution,
called well-balanced schemes.
The reader is referred, for instance,
to~\cite{Jin2001,Xu2002,Audusse2004,Pares2004,Gallardo2007,Berthon2012c}.
More recently, the terminology of fully well-balanced schemes was introduced to
emphasize the ability of some schemes to capture
all linear or nonlinear steady solutions, see e.g.
\cite{Castro2007,Bouchut2010,Berthon2016b,MicBerClaFou2016,Chen2017,Cheng2019,Berthon2024,FraMicNav2024}.
The present work concerns the derivation of a fully well-balanced scheme
for the Euler equations with a gravitational source term.
To that end, we extend techniques which were successfully applied in the context of
fully well-balanced schemes for the shallow water equations,
all the while making sure that the resulting scheme is entropy-stable.

To address the three issues \eqref{i} -- \eqref{ii} -- \eqref{iii},
the paper is organized as follows.
In the next section, we introduce the Euler system under consideration,
and we recall the entropy inequalities satisfied by the weak solutions.
We also present the nonlinear steady solutions with a non-vanishing velocity.
\cref{sec:WBscheme} is then devoted to the derivation of an approximate Riemann solver
which must mimic, in a sense to be prescribed,
the exact Riemann solution associated to the Euler equations with gravitational source terms.
This approximate Riemann solver is defined on the one hand by the integral consistency
condition stated by Harten, Lax and van Leer in~\cite{HarLaxLee1983}.
On the other hand, the approximate Riemann solver must be at rest as soon as a steady solution is considered. This condition is a natural way to suitably define the
discretization of the source term \cite{Berthon2016b,MicBerClaFou2016,Ghitti2020,Berthon2021}.
Next, we impose the entropy integral consistency condition, stated in
\cite{HarLaxLee1983}, to recover the expected discrete entropy
inequalities. Equipped with the approximate Riemann solver, we exhibit
the associated Godunov-type scheme, and we show that the required
properties are satisfied: admissibility of the approximate solutions,
entropy stability, and fully well-balanced property.
In \cref{sec:HigherOrder}, we adopt a recent technique to increase the
order of accuracy of a numerical scheme while conserving the fully well-balanced
property \cite{BerBulFouMbaMic2022}. In \cref{sec:NumRes}, we exhibit
several numerical experiments to illustrate the relevance
of the numerical procedure designed here.
A short conclusion is given in the last section.

\section{The model equations}
\label{sec:model}

This section introduces, in \cref{sec:Euler_gravity},
the Euler equations under a gravitational potential.
Its steady solutions are then derived,
firstly in \cref{sec:moving_equilibria} with a nonzero velocity
and secondly in \cref{sec:hydrostatic_equilibria} with a vanishing velocity,
as the limit of moving steady solutions when the velocity tends to zero.

\subsection{The Euler equations with gravity}
\label{sec:Euler_gravity}

We consider the compressible Euler equations,
equipped with a gravitational source term in a one-dimensional setting,
governed by
\begin{equation}
    \label{eq:EulerG}
    \begin{dcases}
        \pd_t \rho  + \pd_x (\rho u) = 0 ,                             \\
        \pd_t (\rho u)  + \pd_x (\rho u^2 + p) = - \rho \pd_x \varphi, \\
        \pd_t E  + \pd_x ((E + p) u) = - \rho u \pd_x \varphi,
    \end{dcases}
\end{equation}
where $\rho > 0$ denotes the density,
$u \in \mathbb{R}$ the velocity
and $E > 0$ the total energy density.
Moreover, $\varphi : \mathbb{R} \to \mathbb{R}$ is a given time-independent
continuous gravitational potential.
The system is closed by the ideal gas equation of state,
which defines the pressure $p > 0$ by
\begin{equation}
    \label{eq:pressure_law}
    p = (\gamma - 1) \left( E - \frac 1 2 \rho u^2 \right),
\end{equation}
where $\gamma > 1$ denotes the adiabatic index.
To shorten notation, we rewrite system~\eqref{eq:EulerG}
in the following compact form
\begin{equation}
    \label{eq:Euler_compact_form}
    \pd_t W + \pd_x(F(W)) = S(W),
\end{equation}
where the state vector $W$, flux function $F$ and source term $S$ are respectively defined by
\begin{equation*}
    W =
    \begin{pmatrix}
        \rho \\ \rho u \\ E
    \end{pmatrix}, \quad
    F(W) =
    \begin{pmatrix}
        \rho u \\ \rho u^2 + p \\ u (E + p)
    \end{pmatrix}, \text{\quad and \quad}
    S(W) =
    \begin{pmatrix}
        0 \\ - \rho \pd_x \varphi \\ - \rho u \pd_x \varphi
    \end{pmatrix}.
\end{equation*}
In \eqref{eq:EulerG}, we require the density $\rho$ and pressure $p$ to be positive.
Thus, the state vector $W$ belongs to the set of admissible states $\Omega$,
defined by
\begin{equation}
    \label{eq:admissible_states}
    \Omega = \left\{
    W \in \mathbb{R}^3 \text{\; such that \;}
    \rho > 0 \text{\ and \ } p > 0
    \right\}.
\end{equation}

To reflect the influence of the gravitational source term on the wave structure of the Euler equations, we augment system \eqref{eq:EulerG} by the additional equation $\pd_t \varphi = 0$.
The resulting system is hyperbolic with the eigenvalues
\begin{equation*}
    \lambda_{\pm} = u \pm c, \quad
    \lambda_u = u, \quad
    \lambda_0 = 0,
\end{equation*}
where $c$ denotes the speed of sound,
given for the ideal gas law by
\begin{equation*}
    c =
    \sqrt{\frac{\gamma p}{\rho}}.
\end{equation*}
This wave structure is represented in the left panel of \cref{fig:Riemann_solver}.
Note that, compared to the homogeneous system,
which is obtained by setting the gravitational source terms in \eqref{eq:EulerG} to zero,
there is a zero eigenvalue associated to
the gravitational potential.
This leads to a non-ordered wave structure,
since $\lambda_{\pm}$ and $\lambda_u$ can be positive or negative, depending on the flow.
This has a direct consequence on the construction of the Riemann solver,
see e.g.~\cite{Desveaux2016} where a relaxation model is
developed to circumvent the appearance of a zero wave,
or \cite{Vides2014} where six different cases of
wave order configurations had to be taken into account.

Moreover, admissible entropy weak solutions
to system \eqref{eq:EulerG}
satisfy the entropy inequality
\begin{equation}
    \label{eq:entropy_ineq}
    \pd_t (\rho \eta(s))
    +
    \pd_x (\rho u \eta(s)) \leq 0,
\end{equation}
where the specific entropy $s$ is defined by
\begin{equation}
    \label{eq:entropy}
    s = - \log \left( \frac p {\rho^\gamma} \right).
\end{equation}
After \cite{Tadmor2003,Desveaux2016},
the entropy function $W\mapsto \rho\eta(s)$ is convex as soon as
$\eta$ is any smooth function satisfying
\begin{equation*}
    \eta^\prime(s) \geq 0
    \text{\quad and \quad}
    \gamma \eta^{\prime \prime}(s) + \eta^\prime(s) \geq 0.
\end{equation*}
However, after \cite{Tadmor2003,Feireisl2021}, the
entropy stability of \eqref{eq:EulerG} is ensured by merely adopting
entropies defined by smooth convex functions $\eta$ of the state variables such that
\begin{equation}
    \label{eq:entropy_cond}
    \eta^\prime(s) \geq 0
    \text{\quad and \quad}
    \eta^{\prime \prime}(s) \geq 0.
\end{equation}
Also, note that the gravitational source term does not interfere with
the definition of the entropy given above.
However, the presence of gravity in system~\eqref{eq:EulerG}
gives rise to non-trivial steady solutions,
which play an important role in many applications.
Therefore, \cref{sec:hydrostatic_equilibria,sec:moving_equilibria}
are devoted to a derivation of
the steady solutions of system \eqref{eq:EulerG}.

\subsection{Moving equilibria}
\label{sec:moving_equilibria}

Time-invariant solutions of \eqref{eq:EulerG} with a non-zero velocity $u \neq 0$,
the so-called moving equilibria, are governed by the following system
\begin{subequations}\label{eq:EulerG_Stationary}
    \begin{eqnarray}
        \label{eq:EulerG_Stationary.rho}
        && \pd_x (\rho u) = 0 , \\
        \label{eq:EulerG_Stationary.mom}
        && \pd_x (\rho u^2 + p) = - \rho \, \pd_x \varphi, \\
        \label{eq:EulerG_Stationary.energy}
        && \pd_x ((E + p) u) = - \rho u \, \pd_x \varphi.
    \end{eqnarray}
\end{subequations}
Assuming smooth enough steady solutions fulfilling \eqref{eq:EulerG_Stationary}, we obtain from \eqref{eq:EulerG_Stationary.rho} a constant momentum $q_0 = \rho u \neq 0$.
Substituting and dividing by $q_0$
in \eqref{eq:EulerG_Stationary.energy}, we find that
\begin{equation}
    \label{eq:def_constant_enthalpy}
    \frac{E + p}{\rho} + \varphi = H_0,
\end{equation}
where $H_0$ denotes a constant total enthalpy.

Substituting $\varphi$ from \eqref{eq:def_constant_enthalpy}
into \eqref{eq:EulerG_Stationary.mom}, we obtain
\begin{equation}
    \label{eq:EulerG_Stationary.mom2}
    \pd_x \! \left( \frac {q_0^2} \rho + p \right)
    =
    \rho \, \pd_x \! \left( \frac{E + p}{\rho} \right)
\end{equation}
Since $p$ is given by \eqref{eq:pressure_law},
we note that
\begin{equation}
    \label{eq:enthalpy_wrt_rho_q_p}
    E + p
    =
    \frac p {\gamma - 1} + \frac 1 2 \frac{q_0^2}{\rho} + p
    =
    \frac{\gamma p}{\gamma - 1} + \frac 1 2 \frac{q_0^2}{\rho}.
\end{equation}
Using the above relation in
\eqref{eq:EulerG_Stationary.mom2} and rearranging terms yields
\begin{equation*}
    - \frac 1 p \pd_x p
    +
    \gamma \frac 1 \rho \pd_x \rho = 0,
    \text{\qquad i.e., \qquad}
    \pd_x \big( - \log p + \gamma \log \rho \big) = 0.
\end{equation*}%
This immediately leads to
\begin{equation*}
    \pd_x \left( - \log \left( \frac p {\rho^\gamma} \right) \right) = 0,
\end{equation*}
and we recognize the expression of the specific entropy $s$ defined in \eqref{eq:entropy}, and thus $\pd_x s = 0$ holds.

Hence, smooth moving steady solutions obeying the ideal gas law, with $u \neq 0$,
are necessarily isentropic and are characterized by the constant triplet $(q_0,H_0,s_0)$ given by
\begin{equation}
    \label{eq:equilibrium}
    \rho u =: q_0,
    \qquad
    \frac{E + p}{\rho} + \varphi =: H_0,
    \qquad
    - \log \left( \frac p {\rho^\gamma} \right) =: s_0,
\end{equation}
where $q_0$ is a constant nonzero momentum,
$H_0$ is a constant total enthalpy
(including the contribution of the potential energy
to the total energy),
and $s_0$ a constant specific entropy.

\subsection{Hydrostatic equilibria}
\label{sec:hydrostatic_equilibria}

In this section, we shortly recap some properties of
hydrostatic steady solutions for the Euler equations with gravity,
which are equilibria at rest with $u= 0$.
Considering time-invariant solutions and setting $u=0$ in \eqref{eq:EulerG}
leads to steady solutions characterized by the hydrostatic equation
\begin{equation}
    \label{eq:EQ_hydro}
    \pd_x p = - \rho \pd_x \varphi.
\end{equation}
Relation \eqref{eq:EQ_hydro} is the only link between the two unknowns, $p$ and $\rho$.
Thus, hydrostatic equilibria are underdetermined and,
to obtain unique solutions of \eqref{eq:EQ_hydro},
additional information about the properties of the equilibrium
have to be provided.
Intensively studied examples are so-called isothermal or isentropic hydrostatic equilibria with constant temperature or constant entropy, respectively.
The reader is referred to e.g. \cite{Kaepelli2014,Desveaux2016,GrosheintzLaval2019,Thomann2019} for a non-exhaustive list of numerical methods exact on hydrostatic equilibria.

Staying in the framework of the moving steady solutions governed by \eqref{eq:equilibrium},
and taking a vanishing velocity $u\to0$,
the resulting hydrostatic equilibrium is necessarily isentropic, and is governed by
\begin{equation}
    \label{eq:isentropic_hydrostatic_equilibrium}
    \rho u = 0
    \text{, \quad}
    \frac{\kappa \gamma}{\gamma - 1} \rho^{\gamma - 1} + \varphi = \hat{H}_0,
    \quad - \log \left( \frac p {\rho^\gamma} \right) =: s_0,
\end{equation}
where $\kappa = \exp(-s_0)$ is constant.
One can easily check that \eqref{eq:isentropic_hydrostatic_equilibrium}
indeed satisfies the hydrostatic equation \eqref{eq:EQ_hydro}.
Moreover, the usual isentropic hydrostatic steady solution
is recovered as a solution of \eqref{eq:isentropic_hydrostatic_equilibrium}, given by
\begin{equation}
    \label{eq:isentropic_hydrostatic_equilibrium_sol}
    u = 0, \quad p = \kappa \rho^\gamma, \quad \rho = \left(\frac{\gamma - 1}{\gamma \kappa} (\hat{H}_0 - \varphi)\right)^{\frac{1}{\gamma - 1}}.
\end{equation}

In the following, the objective is to derive a Godunov-type finite volume (FV) scheme based on an approximate Riemann solver (\ARS), which is able to preserve moving equilibria and isentropic hydrostatic equilibria up to machine precision,
generates admissible solutions in the sense of \eqref{eq:admissible_states},
and fulfills all entropy inequalities~\eqref{eq:entropy_ineq} with \eqref{eq:entropy_cond}.

\section{The numerical scheme}
\label{sec:WBscheme}

We start by briefly recalling the principle of
Go\-du\-nov-type FV schemes
based on approximate Riemann solvers.
The space domain is discretized, for $i \in \mathbb{Z}$, in uniform cells
$C_i = (x_\imhinline, x_\iphinline)$, of center $x_i$ and of size $\dx$.
The time domain is discretized with time steps of varying size $\dt$.
We seek a way to evolve in time an approximation of the solution $W$
to the Euler equations~\eqref{eq:Euler_compact_form}, defined by
\begin{equation}
    \label{eq:cell_average_W}
    W_i^n \simeq
    \frac 1 \dx
    \int_{x_\imh}^{x_\iph} W(x,t^n)\, dx.
\end{equation}

To that end,
we develop an approximate Riemann solver (\ARS).
We first state, in \cref{sec:Riemann_solver},
the structure of our \ARS,
as well as the general principle of \ARSs.
Then, basic consistency properties are established
in \cref{sec:Riemann_solver_consistency},
and we define the time evolution of the cell averages in
\cref{sec:Godunov_type_scheme},
up to the final choice of the \ARS.
The remainder of this section is then devoted to
the derivation of the remaining unknowns in the \ARS
to achieve the properties given in \cref{sec:definitions}.

\subsection{The approximate Riemann solver}
\label{sec:Riemann_solver}

Let us define an \ARS, denoted by~$\widetilde W$,
whose role is to approximate the solution of the Riemann problems
occurring at each cell interface $x_\iphinline$ between cells $C_i$ and $C_{i+1}$.
Since this solution is self-similar,
the \ARS~depends on the variable $x/t$,
as well as on the two states $W_i^n$ and $W_{i+1}^n$ left and right of the interface.
For instance, the approximate Riemann solution
at interface $x_\iphinline$ is given by
$\smash{\widetilde W((x - x_\iphinline) / t; W_i^n, W_{i+1}^n)}$.
It also depends on the cell averages of the gravitational potential
located left and right of the interface.
These cell averages are denoted by
$\varphi_i$ and $\varphi_{i+1}$,
and are defined similarly to $W_i^n$,
see \eqref{eq:cell_average_W}.
They should also appear as inputs to the \ARS \smash{$\widetilde W$}.
However, to avoid cluttering notation,
we elect not to write these inputs,
and write that \smash{$\widetilde W$} only depends on
the self-similar variable
and on $W_i^n$ and $W_{i+1}^n$.

In this work, we select an \ARS~made of two intermediate states,
as represented in the right panel of \cref{fig:Riemann_solver}.
Its expression, for given left and right states
$W_L$ and $W_R$, reads
\begin{equation}
    \label{eq:ARS}
    \widetilde{W}\left( \frac x t; W_L, W_R \right) =
    \begin{dcases}
        W_L   & \text{if } x < -\lambda t,     \\
        W_L^* & \text{if } -\lambda t < x < 0, \\
        W_R^* & \text{if } 0 < x < \lambda t,  \\
        W_R   & \text{if } x > \lambda t,      \\
    \end{dcases}
\end{equation}
where the components of the vectors $W_L^*$ and $W_R^*$
are denoted by
\begin{equation*}
    W_L^* =
    \begin{pmatrix}
        \rho_L^* \\ q_L^* \\ E_L^*
    \end{pmatrix}
    \text{\qquad and \qquad}
    W_R^* =
    \begin{pmatrix}
        \rho_R^* \\ q_R^* \\ E_R^*
    \end{pmatrix}.
\end{equation*}

\begin{figure}[tb]
    \centering
    \begin{tikzpicture}[scale=0.8]
        \pgfmathsetmacro{\dx}{-7}
        \pgfmathsetmacro{\eps}{0.2}

        \draw[thick, -latex] (-2.5+\dx,0) -- (2.5+\dx,0) node[right]{$x$};
        \draw[thick, -latex] (0+\dx,-0.25) -- (0+\dx,3.25) node[above]{$t$};

        \draw[] (\dx,0) -- (\dx+0.5,2.5) node[above]{$u$};
        \draw (\dx,0) -- (\dx+1.75,2.5) node[above]{$u+c$};
        \draw (\dx,0) -- (\dx-0.75,2.5) node[above]{$u-c$};

        \draw[thick, -latex] (-2.5,0) -- (2.5,0) node[right]{$x$};
        \draw[thick, -latex] (0,-0.25) -- (0,3.25) node[above]{$t$};
        \draw (0,0) -- (-1.75-\eps, 2.5) node[above]{$-\lambda$};
        \draw (0,0) -- (1.75+\eps, 2.5) node[above]{$+\lambda$};

        \node at (-1.6, 0.9) {$W_L$};
        \node at (-0.6, 2) {$W_L^*$};
        \node at (0.6, 2) {$W_R^*$};
        \node at (1.6, 0.9) {$W_R$};
    \end{tikzpicture}

    \caption{%
        Left panel: One possible wave configuration,
        with $u > 0$,
        of the Euler equations with gravity \eqref{eq:EulerG}.
        Right panel: Wave structure of the approximate Riemann solver.%
    }
    \label{fig:Riemann_solver}
\end{figure}

In the \ARS \eqref{eq:ARS}, the approximate wave speed $\lambda$
is chosen such that the acoustic waves are included
in the wave fan of the \ARS,
i.e., $\lambda \geq |u| +c \geq |\lambda_\pm|$.
In this work, we impose that the waves defining the \ARS
are symmetric around $\lambda_0 = 0$.
However, this is just a simplifying assumption.
Indeed, in principle, two different values of $\lambda$
could be chosen, to account for the asymmetric wave configuration
of the Euler equations with respect to $\lambda_0$,
as long the acoustic waves $\lambda_\pm$ are contained in the \ARS.
In practice, we take
\begin{equation}
    \label{eq:def_lambda}
    \lambda = \Lambda \, \max \Big(|u_L| + c(W_L), |u_R| + c(W_R) \Big).
\end{equation}
The parameter $\Lambda \geq 1$ is used
to scale the wave speeds,
and will be given in the numerical experiments.
Indeed, it ensures that the approximate wave speeds in
the approximate Riemann solver are faster than
the true Euler wave speeds.
This is a requirement for stability,
see e.g.~\cite{HarLaxLee1983}.
Therefore, the larger the value of $\Lambda$,
the more diffusion is introduced into the resulting numerical scheme.
This is reflected in the approximate solution by
smeared out shock and contact waves.
Thus, $\Lambda$ is ideally chosen heuristically close to $1$,
to balance between stability and diffusivity of the numerical scheme.

The remainder of this section is devoted to
the construction of the two intermediate states $W_L^*$ and $W_R^*$.
They will be derived such that the resulting numerical scheme satisfies
several key properties,
namely consistency, entropy stability
(i.e. fulfilling a discrete equivalent of \eqref{eq:entropy_ineq}),
positivity preservation (i.e. $W_i^n \in \Omega$),
and well-balancedness.
We start with consistency in
\cref{sec:Riemann_solver_consistency}.

\subsection{Consistency conditions}
\label{sec:Riemann_solver_consistency}

Let us first derive conditions on the intermediate
states~$W_L^*$ and $W_R^*$ to ensure
the consistency of the FV scheme with the solution of the Euler system.
According to~\cite{HarLaxLee1983}, the scheme is consistent
as soon as the average over a cell of the \ARS~\eqref{eq:ARS}
is equal to that of the exact, self-similar solution to the
Riemann problem, denoted by $W_\mathcal{R}$.
Therefore, the following equality must hold
\begin{equation}
    \label{eq:integral_consistency}
    \frac 1 \dx \int_{-\dx / 2}^{\dx / 2}
    \widetilde{W} \left( \frac x \dt; W_L, W_R \right) dx
    =
    \frac 1 \dx \int_{-\dx / 2}^{\dx / 2}
    W_\mathcal{R} \left( \frac x \dt; W_L, W_R \right) dx.
\end{equation}
Straightforward computations,
detailed in \cite{MicBerClaFou2016},
then show that the intermediate states must satisfy the following relations:
\begin{equation}
    \label{eq:raw_consistency_condition}
    \left\{
    \begin{aligned}
        \rho_L^* + \rho_R^* & = 2 \rho_\HLL,      \\
        q_L^* + q_R^*       & =
        2 q_\HLL +
        \frac 1 {\lambda \dt}
        \int_{-\dx / 2}^{\dx / 2}
        \int_0^{\dt}
        (-\rho \pd_x \varphi)_\mathcal{R}
        \left( \frac x t; W_L, W_R \right) dt dx, \\
        E_L^* + E_R^*       & =
        2 E_\HLL +
        \frac 1 {\lambda \dt}
        \int_{-\dx / 2}^{\dx / 2}
        \int_0^{\dt}
        (-q \pd_x \varphi)_\mathcal{R}
        \left( \frac x t; W_L, W_R \right) dt dx, \\
    \end{aligned}
    \right.
\end{equation}
where $(-\rho \pd_x \varphi)_\mathcal{R}$ and
$(-q \pd_x \varphi)_\mathcal{R}$ respectively represent the
second and third components of $S(W_\mathcal{R})$.
Moreover, we have introduced in \eqref{eq:raw_consistency_condition}
the intermediate state~$W_\HLL$ of the well-known HLL Riemann solver
from Harten, Lax and van Leer~\cite{HarLaxLee1983}.
It is given by
\begin{equation}
    \label{eq:def_w_hll}
    W_\HLL =
    \begin{pmatrix}
        \rho_\HLL \\ q_\HLL \\ E_\HLL
    \end{pmatrix} =
    \frac{W_L + W_R} 2 - \frac{F(W_R) - F(W_L)} {2 \lambda}.
\end{equation}

Next, we introduce two new quantities, $S^q$ and $S^E$,
which are approximations of the cell averages of the integrals
of the second and third components of the source term, i.e.,
\begin{equation}
    \label{eq:consistent_cell_averages_of_source_term}
    \begin{aligned}
        S^q & \approx
        \frac 1 {\dt}
        \frac 1 {\dx}
        \int_{-\dx / 2}^{\dx / 2}
        \int_0^{\dt}
        (-\rho \pd_x \varphi)_\mathcal{R}
        \left( \frac x t; W_L, W_R \right) dt dx, \\
        S^E & \approx
        \frac 1 {\dt}
        \frac 1 {\dx}
        \int_{-\dx / 2}^{\dx / 2}
        \int_0^{\dt}
        (-q \pd_x \varphi)_\mathcal{R}
        \left( \frac x t; W_L, W_R \right) dt dx, \\
    \end{aligned}
\end{equation}
where $S^q$ and $S^E$ have to be consistent with
$-\rho \pd_x \varphi$ and $-q \pd_x \varphi$,
in a sense that will be prescribed later.
Then, using relations \eqref{eq:consistent_cell_averages_of_source_term}
in the consistency condition~\eqref{eq:raw_consistency_condition}
yields the following consistency relations to be satisfied by the
intermediate states
\begin{equation}
    \label{eq:consistency_condition_with_source_averages}
    \left\{
    \begin{aligned}
        \rho_L^* + \rho_R^* & = 2 \rho_\HLL,  \\
        q_L^* + q_R^*       & =
        2 q_\HLL + \frac {S^q \dx} {\lambda}, \\
        E_L^* + E_R^*       & =
        2 E_\HLL + \frac {S^E \dx} {\lambda}. \\
    \end{aligned}
    \right.
\end{equation}

To simplify computations in the following steps,
let us introduce $\widehat W$ and $\delta W$,
defined such that
\begin{equation}
    \label{eq:definition_w_hat_and_delta_w_wrt_w_star}
    \widehat W = \frac{W_L^* + W_R^*}2
    \text{\qquad and \qquad}
    \delta W = \frac{W_R^* - W_L^*}2.
\end{equation}
We highlight that this corresponds to rewriting
the intermediate states in terms of the unknowns
$\widehat W$ and $\delta W$.
Namely, we have $W_L^* = \widehat W - \delta W$
and $W_R^* = \widehat W + \delta W$.
Equipped with this notation,
\eqref{eq:consistency_condition_with_source_averages}
immediately gives the following expressions for the components of $\widehat{W}$:
\begin{equation}
    \label{eq:expression_w_hat}
    \left\{
    \begin{aligned}
        \widehat \rho & = \rho_\HLL,   \\
        \widehat q    & = q_\HLL
        + \frac {S^q \dx} {2 \lambda}, \\
        \widehat E    & = E_\HLL
        + \frac {S^E \dx} {2 \lambda}. \\
    \end{aligned}
    \right.
\end{equation}
We are thus left with the determination of the three components of $\delta W$
and of the two source term averages $S^q$ and $S^E$,
i.e. five unknowns in total.
However, whatever the expressions of these unknowns,
the scheme \eqref{eq:scheme_from_ARS}
is constructed to satisfy \eqref{eq:integral_consistency},
and is therefore consistent by definition.

The remainder of \cref{sec:WBscheme} is devoted to the derivation of these five unknowns,
ensuring positivity preservation, entropy stability and well-balancedness.
However, first, we explain in \cref{sec:Godunov_type_scheme}
how to use the \ARS \eqref{eq:ARS} to derive the time update
of the numerical scheme,
and precisely define in \cref{sec:definitions}
the properties to be satisfied by the scheme,
and thus by the intermediate states.

\subsection{The Godunov-type scheme finite volume scheme}
\label{sec:Godunov_type_scheme}

We start by defining the cell-wise juxtaposition $W_\Delta$
of the approximate Riemann solver \eqref{eq:ARS}, by
\begin{equation*}
    \forall i \in \mathbb{Z}, \,
    \forall t \in (0, \dt], \,
    \forall x \in (x_i, x_{i+1}), \,
    W_\Delta(x, t^n + t)
    =
    \widetilde W\left(\dfrac{x - x_\iph}{t}; W_{i}^n, W_{i+1}^n\right).
\end{equation*}
This definition is illustrated on
\cref{fig:juxtaposition_of_ARSs}.
The time update of the Godunov-type scheme is then given by
integrating the juxtaposition $W_\Delta$ of Riemann solvers
at time $t^{n+1} = t^n + \dt$ over each cell $C_i$, yielding
\begin{equation}
    \label{eq:wnp1}
    W_i^{n+1}=\frac{1}{\dx}\int_{x_\imh}^{x_\iph}
    W_\Delta(x,t^{n+1})dx.
\end{equation}
For stability purposes,
we impose the following CFL (Courant-Friedrichs-Lewy)
condition on $\dt$, see \cite{HarLaxLee1983}:
\begin{equation}
    \label{eq:CFL_condition}
    \dt \leq \frac 1 2 \frac{\dx}{\max_i \lambda_\iph}.
\end{equation}

\begin{figure}[tb]
    \centering
    \begin{tikzpicture}[scale=0.9]
        \draw[thick, -latex] (-0.75,0) node [left]{$t^n$} -- (10.5,0) node [right]{$x$};
        \draw[dashed] (-0.75,3) node [left]{$t^{n+1}$} -- (10.5,3);

        \draw[thick, -latex] (-0.5,-0.25) -- (-0.5,3.5) node[above]{$t$};

        \draw[densely dotted, thick] (0,-0.25) node[below]{$x_\imh$} -- (0,3) node (x_imh) {};
        \draw[densely dotted, thick] (10,-0.25) node[below]{$x_\iph$} -- (10,3) node (x_iph) {};
        \node at (5,0) {$\times$};
        \node[below] at (5,-0.25) {$x_i$};
        \node (x_i) at (5,3) {};

        \draw (0,0) -- (2.5,3) node[midway, right]{$\quad\lambda_\imh$};
        \draw (10,0) -- (6.5,3) node[midway, left]{$-\lambda_\iph \quad$};

        \node [below] at (1, 2.75) {$W_{\imh, R}^*$};
        \node [below] at (5, 2.75) {$W_i^n$};
        \node [below] at (8.5, 2.75) {$W_{\iph, L}^*$};

        \draw [decoration={brace, amplitude=5pt}, decorate] (x_imh.north) -- (x_i.north)
        node [midway, above, font=\small, yshift=2pt]
        {$\widetilde W\left(\dfrac{x - x_\imh}{\dt}; W_{i-1}^n, W_i^n\right)$};

        \draw [decoration={brace, amplitude=5pt}, decorate] (x_i.north) -- (x_iph.north)
        node [midway, above, font=\small, yshift=2pt]
        {$\widetilde W\left(\dfrac{x - x_\iph}{\dt}; W_{i}^n, W_{i+1}^n\right)$};

        \draw [decoration={brace, amplitude=5pt, raise=30pt}, decorate] (x_imh.north) -- (x_iph.north)
        node [midway, above, font=\small, yshift=32pt]
        {$W_\Delta(x, t^{n} + \dt)$};
    \end{tikzpicture}
    \caption{%
        Juxtaposition of subsequent \ARSs
        defining $W_\Delta$ using definition \eqref{eq:IS_short} of the intermediate states
        $W_{\imhinline, R}^*$ and $W_{\iphinline, L}^*$. %
    }
    \label{fig:juxtaposition_of_ARSs}
\end{figure}

To properly compute the integral in \eqref{eq:wnp1},
we note that $W_L^*$ and~$W_R^*$
are given as functions
of the left and right states $W_L$ and~$W_R$,
as well as of the corresponding left and right
cell averages of the gravitational potential, $\varphi_L$ and~$\varphi_R$.
Since the structure of the \ARS~\eqref{eq:ARS}
comprises a stationary wave with velocity~$0$,
the scheme, see e.g. \cite{MicBerClaFou2016}, is given by
\begin{equation}
    \label{eq:scheme_from_ARS}
    W_i^{n+1} = W_i^n + \frac \dt \dx \Big[
        \lambda_{\iph} \Big(W_L^* \left( W_i^{n}, W_{i+1}^n, \varphi_i, \varphi_{i+1} \right) - W_i^n\Big)
        +
        \lambda_{\imh} \Big(W_R^* \left( W_{i-1}^{n}, W_i^n, \varphi_{i-1}, \varphi_i \right) - W_i^n\Big)
        \Big].
\end{equation}

We emphasize that the scheme \eqref{eq:scheme_from_ARS}
can also be rewritten in the standard con\-ser\-va\-tion form
\begin{equation}
    \label{eq:scheme_in_conservation_form}
    W_i^{n+1} = W_i^n
    -
    \frac{\dt}{\dx} \left(
    \mathcal{F}(W_{i}^n, W_{i+1}^n)
    -
    \mathcal{F}(W_{i-1}^n, W_{i}^n)
    \right)
    +
    \frac{\dt}{2} \left(
    \mathcal{S}_\imh^n + \mathcal{S}_\iph^n
    \right),
\end{equation}
where $\mathcal{F}(W_{i}^n, W_{i+1}^n)$
is the numerical flux at interface $x_\iphinline$,
and where $\mathcal{S}_\iphinline^n$
is the numerical source term at interface $x_\iphinline$.
The expressions of $\mathcal{F}(W_{i}^n, W_{i+1}^n)$
and~$\mathcal{S}_\iphinline^n$
are obtained thanks to straightforward computations
presented in e.g.~\cite{HarLaxLee1983}.
Setting
\begin{align}
    \label{eq:IS_short}
    \begin{split}
         & W_{\iphinline, L}^* = W_L^*( W_{i-1}^{n}, W_i^n, \varphi_{i-1}, \varphi_i) \text{ \quad and} \\
         & W_{\iphinline, R}^* = W_R^*( W_{i-1}^{n}, W_i^n, \varphi_{i-1}, \varphi_i),
    \end{split}
\end{align}
we obtain
\begin{equation}
    \label{eq:num_flux}
    \mathcal{F}(W_{i}^n, W_{i+1}^n)
    =
    \frac {F(W_i^n) + F(W_{i+1}^n)} 2
    -
    \frac {\lambda_\iph}{2} \left(W_{\iph, L}^* - W_i^n \right)
    +
    \frac {\lambda_\iph}{2} \left(W_{\iph, R}^* - W_{i+1}^n \right).
\end{equation}
The numerical source term is given by
\begin{equation*}
    \mathcal{S}_\iph^n
    =
    \begin{pmatrix}
        0                         \\
        \left( S^q \right)_\iph^n \\
        \left( S^E \right)_\iph^n
    \end{pmatrix},
\end{equation*}
where \smash{$(S^q)_\iphinline^n$}
and \smash{$(S^E)_\iphinline^n$} denote
the source term approximations
\eqref{eq:consistent_cell_averages_of_source_term}
evaluated at the interface $x_\iphinline$.
The numerical flux defined in \eqref{eq:num_flux}
also depends on the left and right cell averages of the
gravity potential, $\varphi_i$ and $\varphi_{i+1}$.
We recall that this dependency is not explicitly written,
to avoid cluttered notation.

By \eqref{eq:scheme_in_conservation_form},
or equivalently \eqref{eq:scheme_from_ARS},
we have defined the time update of the consistent
Godunov-type finite volume scheme,
up to the definition of the
intermediate states $W_L^*$ and $W_R^*$
and source term approximations $S^q$ and $S^E$.
Before giving their derivation,
we detail in \cref{sec:definitions} the key properties the scheme must satisfy.

\subsection{Definitions and conditions to be satisfied by the \ARS}
\label{sec:definitions}

The objective of this section is to state important definitions
and give theoretical results regarding the well-balanced property, positivity and entropy stability.
They are fundamental for the derivation of the unknowns in the \ARS,
which will be carried out in \cref{sec:determination_of_delta_w,sec:source_term_approximation}.

\subsubsection{Well-balancedness}

We start with the definition of a well-balanced scheme.
We emphasize that we are interested in developing a so-called
\emph{fully well-balanced} scheme,
which exactly preserves the moving steady solutions
described in \cref{sec:moving_equilibria} and isentropic hydrostatic equilibria as consequence.

To properly introduce our definition of well-balancedness,
we first need to introduce the notion of
Interface Steady Solution (ISS).
Indeed, from \eqref{eq:equilibrium}, we know that a solution is steady
if the momentum $q$, the total enthalpy $H$
and the specific entropy $s$ are constant.
This leads us to the following definition of an ISS.

\begin{definition}[Interface Steady Solution (ISS)]
    \label{def:ISS}
    A pair $(W_L, W_R)$ of admissible states is said to
    be an Interface Steady Solution (ISS) if, and only if,
    \begin{equation*}
        q_L = q_R
        \text{, \quad}
        H(W_L, \varphi_L) = H(W_R, \varphi_R)
        \text{, \quad and \quad}
        s(W_L) = s(W_R).
    \end{equation*}
\end{definition}

Equipped with the definition of the ISS,
we can now properly define a well-balanced scheme.

\begin{definition}[Well-balancedness]
    \label{def:well-balancedness}
    Scheme \eqref{eq:scheme_from_ARS} is said to be well-balanced if
    \begin{equation*}
        \forall i \in \mathbb{Z}, \; (W_i^n, W_{i+1}^n) \text{ is an ISS}
        \qquad \implies \qquad
        \forall i \in \mathbb{Z}, \; W_i^{n+1} = W_i^n.
    \end{equation*}
\end{definition}

\begin{lemma}
    \label{cor:WB_consequence}
    Let $(W_L, W_R)$ be an ISS.
    A sufficient condition for well-ba\-lan\-ced\-ness is that
    $W_L^* = W_L$ and $W_R^* = W_R$.
\end{lemma}

\begin{proof}
    Assume that
    $(W_i^n, W_{i+1}^n)$ and $(W_{i-1}^n, W_i^n)$ are ISSs.
    We need to show that $W_i^{n+1} = W_i^n$
    in order for the scheme to preserve the steady solution.
    If the intermediate states satisfy
    $W_L^* = W_L$ and $W_R^* = W_R$, then
    $W_L^*( W_i^{n}, W_{i+1}^n, \varphi_i, \varphi_{i+1}) = W_i^n$
    and
    $W_R^* ( W_{i-1}^{n}, W_i^n, \varphi_{i-1}, \varphi_i ) = W_i^n$.
    Plugging these relations into
    \eqref{eq:scheme_from_ARS} immediately yields
    $W_i^{n+1} = W_i^n$, which concludes the proof.
\end{proof}

According to \cref{cor:WB_consequence},
and recalling the definitions
\eqref{eq:definition_w_hat_and_delta_w_wrt_w_star}
of $\widehat W$ and $\delta W$,
it is sufficient that the following conditions hold
\begin{equation}
    \label{eq:WB_consequence_on_hat_and_delta_w}
    (W_L, W_R) \text{ is an ISS }
    \quad \implies \quad
    \widehat W = \frac{W_L + W_R} 2
    \text{ \; and \; }
    \delta W = \frac{W_R - W_L} 2.
\end{equation}
Put in other words, $(W_L^*, W_R^*)$ should also be an ISS, i.e.,
\begin{equation}
    \label{eq:WB_consequence_on_intermediate_states}
    q_L^* = q_R^*
    \text{, \qquad}
    H(W_L^*, \varphi_L) = H(W_R^*, \varphi_R)
    \text{\qquad and \qquad}
    s(W_L^*) = s(W_R^*).
\end{equation}
Conditions~\eqref{eq:WB_consequence_on_intermediate_states}
lead for instance to a single momentum component for both intermediate states,
which we denote by $q^* \coloneqq q_L^* = q_R^*$.
Moreover, they define a single intermediate entropy
\begin{equation}
    \label{eq:intermediate_entropy}
    s^* \coloneqq s(W_L^*) = s(W_R^*).
\end{equation}
Conditions
\eqref{eq:WB_consequence_on_hat_and_delta_w} and
\eqref{eq:WB_consequence_on_intermediate_states}
will be instrumental when determining the expressions
of $\delta W$ and  $S^q$,~$S^E$
in \cref{sec:determination_of_delta_w,sec:source_term_approximation}, respectively.

\subsubsection{Entropy stability}

The second crucial property that we wish to satisfy
is entropy stability, defined below.

\begin{definition}[Entropy stability]
    \label{def:entropy_stability}
    The scheme \eqref{eq:scheme_from_ARS} is said to be entropy-stable if,
    for all smooth functions $\eta$ satisfying \eqref{eq:entropy_cond}
    and for all $i\in\mathbb{Z}$,
    the states $W_i^n$ and $W_i^{n+1}$ fulfill the discrete entropy inequality
    \begin{equation}\label{discreteentropineq}
        \rho_i^{n+1}\eta(s_i^{n+1}) \leq \rho_i^n \eta(s_i^n)
        - \frac{\dt}{\dx} \left((\rho \eta(s) u)_{i+1/2}^n -
        (\rho \eta(s) u)_{i-1/2}^n\right).
    \end{equation}
\end{definition}

In order to reach this important stability property and since we are
considering Godunov-type schemes, we here adopt the well-known
integral entropy consistency stated in \cite{HarLaxLee1983}. As a consequence,
the scheme \eqref{eq:scheme_from_ARS} will be entropy-preserving by
satisfying a discrete entropy inequality \eqref{discreteentropineq} as
soon as the \ARS, given by $\smash{\widetilde{W}}$ in the form \eqref{eq:ARS},
verifies
\begin{equation}
    \label{eq:entropy_integral_consistency}
    \frac{1}{\dx}\int_{-\dx/2}^{\dx/2}
    \widetilde{(\rho\eta(s))}  \left( \frac{x}{\dt},W_L,W_R\right) dx
    \le
    \frac{1}{2}\big( \rho_L    \eta(s_L) + \rho_R\eta(s_R) \big)
    -\frac{\dt}{\dx}
    \big( \rho_R\eta(s_R)u_R + \rho_L\eta(s_L)u_L \big),
\end{equation}
where, with abuse of notation, we have defined $s_L=s(W_L)$
and $s_R=s(W_R)$.
According to \eqref{eq:ARS}, we readily show that the integral
in \eqref{eq:entropy_integral_consistency} satisfies
\begin{equation*}
    \begin{aligned}
        \frac{1}{\dx}
        \int_{-\dx/2}^{\dx/2}
        \widetilde{(\rho\eta(s))}
        \left( \frac{x}{\dt},W_L,W_R\right) dx
         & =
        \left(\frac{1}{2} - \lambda\frac{\dt}{\dx}\right)
        \rho_L\eta(s_L)
        +\lambda\frac{\dt}{\dx} \rho^*_L\eta(s(W^*_L))    \\
         & +\lambda\frac{\dt}{\dx} \rho^*_R\eta(s(W^*_R))
        + \left(\frac{1}{2} - \lambda\frac{\dt}{\dx}\right)
        \rho_R\eta(s_R),
    \end{aligned}
\end{equation*}
so that, by involving $s(W^*_L) = s(W^*_R) \eqqcolon s^*$
imposed in \eqref{eq:intermediate_entropy},
\eqref{eq:entropy_integral_consistency} is immediately rewritten as
\begin{equation*}
    \frac{1}{2}\big( \rho^*_L+\rho^*_R\big) \eta(s^*) \le
    \frac{1}{2}\big( \rho_L\eta(s_L) + \rho_R\eta(s_R) \big)
    -\frac{1}{2\lambda}
    \big( \rho_R\eta(s_R)u_R + \rho_L\eta(s_L)u_L \big).
\end{equation*}
For the sake of simplicity in forthcoming notation, in addition
to \eqref{eq:def_w_hll}, we define
\begin{align}
    \label{defsHLL}
    (\rho s)_\HLL                 & = \frac{1}{2}
    \big(\rho_Ls_L + \rho_Rs_R\big)
    -\frac{1}{2\lambda}\big(\rho_Rs_Ru_R -\rho_Ls_Lu_L\big), \\
    \label{defetasHLL}
    \big(\rho \eta(s)\big)_{\HLL} & = \frac{1}{2}
    \big(\rho_L\eta(s_L) + \rho_R\eta(s_R)\big)
    -\frac{1}{2\lambda}\big(\rho_R\eta(s_R)u_R -\rho_L\eta(s_L)u_L\big).
\end{align}
As a consequence of this notation,
and because \eqref{eq:consistency_condition_with_source_averages}
implies that $\rho_L^* + \rho_R^* = 2 \rho_\HLL$,
the integral entropy consistency
condition \eqref{eq:entropy_integral_consistency} to be satisfied by the \ARS now reads
\begin{equation}\label{ineqentrop2}
    \rho_{\HLL}\eta(s^*)\le \big(\rho \eta(s)\big)_{\HLL}.
\end{equation}
In order to select a suitable definition of $s^*$, we now state one of
our main results.
\begin{theorem}\label{th:entrop}
    Let
    $W_L \in \Omega$ and $W_R \in \Omega$ be two given states,
    and let us adopt a definition of $\lambda$ such that $\rho_{\HLL}>0$.
    Then, for all convex smooth functions $\eta$, we have
    \begin{equation}\label{ineqentropHLLeta}
        \eta\left(\frac{(\rho s)_\HLL}{\rho_{\HLL}}\right)\le
        \frac{\big(\rho \eta(s)\big)_{\HLL}}{\rho_{\HLL}}.
    \end{equation}
\end{theorem}
\begin{proof}
    The establishment of this important result relies on the introduction
    of a standard HLL scheme for an auxiliary system given by
    \begin{subequations}
        \label{sysaux}
        \begin{eqnarray}
            && \pd_t \rho + \pd_x (\rho u) = 0 ,\\
            && \pd_t (\rho u) + \pd_x (\rho u^2 + p) = 0, \\
            \label{eqrhosaux}
            && \pd_t(\rho s) + \pd_x (\rho s u) = 0,
        \end{eqnarray}
    \end{subequations}
    where the pressure law is given by $s=-\log(p/\rho^\gamma)$ as in \eqref{eq:entropy}. We remark
    that the weak solutions of the above system satisfy, for all smooth
    functions $g$,
    \begin{equation}\label{eqadditional}
        \pd_t(\rho g(s)) + \pd_x (\rho g(s) u) = 0.
    \end{equation}
    We introduce ${U}_\mathcal{R}(x/t;W_L,W_R)$ the exact Riemann solution of
    \eqref{sysaux}, and we consider an \ARS, denoted by $\widetilde{U}(x/t;W_L,W_R)$,
    made of a single intermediate state
    $(\xoverline{\rho}_\HLL,\xoverline{\rho u}_\HLL,\xoverline{\rho
            s}_\HLL)$. Once again, by considering the integral consistency
    condition, here given by
    \begin{equation*}
        \frac{1}{\dx} \int_{-\dx/2}^{\dx/2}
        \widetilde{U}\left(\frac{x}{\dt};W_L,W_R\right) dx
        =
        \frac{1}{\dx} \int_{-\dx/2}^{\dx/2}
        {U}_\mathcal{R}\left(\frac{x}{\dt};W_L,W_R\right) dx,
    \end{equation*}
    we obtain after straightforward computations
    \begin{equation*}
        \begin{aligned}
             & \xoverline{\rho}_\HLL = \frac{1}{2}\left( \rho_L+\rho_R\right)
            - \frac{1}{2\lambda}\left( \rho_R u_R - \rho_L u_L \right),               \\
             & \xoverline{\rho s}_\HLL = \frac{1}{2}\left( \rho_Ls_L+\rho_Rs_R\right)
            - \frac{1}{2\lambda}\left( \rho_R s_R u_R - \rho_L s_L u_L
            \right).
        \end{aligned}
    \end{equation*}
    According to notation \eqref{eq:def_w_hll} and \eqref{defsHLL}, we
    remark that
    \begin{equation}\label{eqnotationHLL}
        \xoverline{\rho}_\HLL = \rho_\HLL
        \text{\qquad and \qquad}
        \xoverline{\rho s}_\HLL = (\rho s)_\HLL.
    \end{equation}
    Next, we integrate
    the equation \eqref{eqrhosaux}, governing $\rho s$, over a space-time slab to get
    \begin{equation*}
        \int_0^{\dt}\int_{-\lambda\dt}^{\lambda\dt}
        \big(
        \pd_t(\rho s) + \pd_x (\rho s u)
        \big)dx\:dt=0.
    \end{equation*}
    Computing the integral, we obtain
    \begin{equation*}
        \xoverline{\rho s}_\HLL = \frac{1}{2\lambda\dt}
        \int_{-\lambda\dt}^{\lambda\dt}
        (\rho s)_\mathcal{R}\left(\frac{x}{\dt};W_L,W_R\right)dx,
    \end{equation*}
    where $(\rho s)_\mathcal{R}$ is the last component of $U_\mathcal{R}$. Then, we have
    \begin{equation*}
        \frac{\xoverline{\rho s}_\HLL}{\xoverline{\rho}_\HLL}
        =
        \frac{1}{2\lambda\dt}
        \int_{-\lambda\dt}^{\lambda\dt}
        s_\mathcal{R}\left(\frac{x}{\dt};W_L,W_R\right)
        \frac{\rho_\mathcal{R}\left(\frac{x}{\dt};W_L,W_R\right)dx}{\xoverline{\rho}_\HLL},
    \end{equation*}
    where we have set $s_\mathcal{R}=(\rho s)_\mathcal{R}/\rho_\mathcal{R}$.
    Now, remark that the measure
    \begin{equation*}
        \frac{\rho_\mathcal{R}\left(\frac{\cdot }{\dt};W_L,W_R\right)dx}
        {2\lambda\dt\xoverline{\rho}_\HLL}
    \end{equation*}
    is a probability measure over $(-\lambda \dt, \lambda \dt)$.
    Hence, for all convex real functions $\eta$,
    Jensen's inequality gives
    \begin{equation}\label{ineqentropaux}
        \eta\left(\frac{\xoverline{\rho s}_\HLL}{\xoverline{\rho}_\HLL}\right)
        \le
        \frac{1}{2\lambda\dt}
        \int_{-\lambda\dt}^{\lambda\dt}
        \eta\left(s_\mathcal{R}\left(\frac{x}{\dt};W_L,W_R\right)\right)
        \frac{\rho_\mathcal{R}\left(\frac{x}{\dt};W_L,W_R\right)dx}{\xoverline{\rho}_\HLL}.
    \end{equation}
    Next, since $U_\mathcal{R}$ defines an exact solution of \eqref{sysaux},
    from \eqref{eqadditional} we get, with $g = \eta$,
    \begin{equation*}
        \int_0^{\dt}\int_{-\lambda\dt}^{\lambda\dt}
        \big(
        \pd_t(\rho \eta(s)) + \pd_x (\rho \eta(s) u)
        \big)dx\:dt=0.
    \end{equation*}
    We thus easily obtain
    \begin{align*}
        \frac{1}{2\lambda\dt}
        \int_{-\lambda\dt}^{\lambda\dt}
        (\rho \eta(s))_\mathcal{R}\left(\frac{x}{\dt};W_L,W_R\right)dx
          =
        \frac{1}{2}\left( \rho_L\eta(s)_L+\rho_R\eta(s)_R\right)
          - \frac{1}{2\lambda}\left( \rho_R \eta(s)_R u_R - \rho_L \eta(s)_L u_L
        \right).
    \end{align*}
    Using notation \eqref{defetasHLL}, we have
    \begin{equation*}
        \frac{1}{2\lambda\dt}
        \int_{-\lambda\dt}^{\lambda\dt}
        (\rho \eta(s))_\mathcal{R}\left(\frac{x}{\dt};W_L,W_R\right)dx
        =(\rho \eta(s))_\HLL
    \end{equation*}
    so that the relation \eqref{ineqentropaux} now reads
    \begin{equation*}
        \eta\left(\frac{\xoverline{\rho s}_\HLL}{\xoverline{\rho}_\HLL}\right)
        \le
        \frac{(\rho\eta(s))_\HLL}{\xoverline{\rho}_\HLL}.
    \end{equation*}
    Using the notation equivalence \eqref{eqnotationHLL}, the estimate
    \eqref{ineqentropHLLeta} is established, and the proof is thus completed.
\end{proof}


Now, using \cref{th:entrop}, we remark that the required
integral entropy consistency condition \eqref{ineqentrop2} is
immediately satisfied as soon as
\begin{equation}
    \label{eq:condition_on_s_star}
    \rho_\HLL s^* = (\rho s)_\HLL,
\end{equation}
which defines the intermediate entropy $s^*$.

\subsubsection{Positivity preservation}
After defining entropy stability,
we are now able to precisely state the
third required property, positivity preservation,
and to give additional conditions on the intermediate states.
Indeed, since $\Omega$ given by \eqref{eq:admissible_states} is a non-convex set,
entropy stability is crucial for a part
of the positivity proof.

\begin{definition}[Positivity preservation]
    \label{def:positivity_preservation}
    Let $W_i^n \in \Omega$ for all $i \in \mathbb{Z}$ and $n \geq 0$.
    The scheme \eqref{eq:scheme_from_ARS} is said to be positivity-preserving if,
    for all $i \in \mathbb{Z}$, $W_i^{n+1} \in \Omega$,
    where $\Omega$ is the set of admissible states defined in \eqref{eq:admissible_states}.
\end{definition}

We first derive a sufficient condition for the positivity of the density.
To that end, recall that $\rho_\HLL>0$ under suitable conditions on $\lambda$.
Indeed, we have
\begin{equation*}
    \rho_\HLL = \frac{\rho_L}{2\lambda} (\lambda+u_L)
    + \frac{\rho_R}{2\lambda} (\lambda-u_R).
\end{equation*}
As a consequence, with $\rho_L>0$ and $\rho_R>0$, the required positivity of $\rho_\HLL$ is
satisfied as soon as $\lambda>\max(|u_L|,|u_R|)$,
which is the case here since $\lambda$ is given
by~\eqref{eq:def_lambda}.
Next, a way to ensure positivity preservation on $\rho_i^{n+1}$
is to impose that $\rho_L^* > 0$ and $\rho_R^* > 0$.
Indeed, in this case, all densities in the Riemann fan are positive,
which implies the positivity of the density at time $t^{n+1}$.
Finally, recall that $\rho_L^* = \rho_\HLL - \delta\rho$
and $\rho_R^* = \rho_\HLL + \delta\rho$.
Therefore, a sufficient condition for the positivity of density is
\begin{equation}
    \label{eq:positivity_condition_on_rho}
    |\delta \rho| < \rho_\HLL.
\end{equation}

Now, with $\rho^*_L>0$ and $\rho^*_R>0$, under assumption \eqref{eq:positivity_condition_on_rho},
we may also exhibit a natural
condition to get a positive updated pressure $p(W_i^{n+1})$.
In fact, we use the discrete entropy inequality to address this issue.
\begin{lemma}[Positivity of the pressure]
    \label{lem:positivity_pressure}
    Let $W_L \in \Omega$ and $W_R \in \Omega$ two given states,
    and let us adopt a definition of $\lambda$ such that $\rho_{\HLL}>0$.
    Moreover, let $\rho_L^* >0 $ and $\rho_R^* > 0$,
    and consider $s^*$ given by \eqref{eq:condition_on_s_star}.
    Then, for all $i \in \mathbb{Z}$, $p_i^{n+1} > 0$.
\end{lemma}
\begin{proof}
    Since $W\mapsto \rho\eta(s)$ is convex, by applying Jensen's inequality
    to the integral form of $W_i^{n+1}$ given by \eqref{eq:wnp1}, we
    obtain
    \begin{equation*}
        \rho_i^{n+1}\eta(s(W_i^{n+1}))\le
        \frac{1}{\dx}\int_{x_\imh}^{x_\iph}
        \rho_\Delta(x,t^{n+1})\eta(s_\Delta(x,t^{n+1}))dx,
    \end{equation*}
    where we have set $s_\Delta = s(W_\Delta)$.
    Using the definition \eqref{eq:entropy} of $s$, we have
    \begin{equation}
        \label{eq:eta_of_minus_log}
        \eta \left( -\log \left( \frac{p(W_i^{n+1})}{\rho_i^{n+1}} \right) \right) \le
        \int_{x_\imh}^{x_\iph}
        \eta \left( -\log \left(
            \frac{p(W_\Delta(x,t^{n+1}))}{(\rho_\Delta(x,t^{n+1}))^\gamma}
            \right) \right)
        \frac{\rho_\Delta(x,t^{n+1}) dx}{\dx \rho_i^{n+1}}.
    \end{equation}
    Moreover,
    we remark that, since $\rho_i^{n+1} > 0$, the measure defined by
    \begin{equation*}
        \frac{\rho_\Delta(\cdot,t^{n+1})dx}{\dx \rho_i^{n+1}}
        =
        \frac{\rho_\Delta(\cdot,t^{n+1})dx}
        {\int_{x_\imh}^{x_\iph}
            \rho_\Delta(x,t^{n+1})dx}
    \end{equation*}
    is a probability measure over $(x_\imhinline,x_\iphinline)$.
    Now, let us introduce the function
    \begin{equation*}
        \begin{aligned}
            \Xi : (0, +\infty) & \to \mathbb{R}          \\
            z                  & \mapsto \eta(-\log(z)).
        \end{aligned}
    \end{equation*}
    Equipped with the properties \eqref{eq:entropy_cond} of $\eta$,
    it is easy to show that $\Xi$ is
    a monotonically decreasing and convex function.
    Therefore, it is invertible, and its inverse $\Xi^{-1}$ is
    monotonically decreasing and concave.
    Applying $\Xi^{-1}$ to \eqref{eq:eta_of_minus_log} yields
    \begin{equation*}
        \frac{p(W_i^{n+1})}{\rho_i^{n+1}} \ge
        \Xi^{-1} \left(
        \int_{x_\imh}^{x_\iph}
        \Xi \left(
            \frac{p(W_\Delta(x,t^{n+1}))}{(\rho_\Delta(x,t^{n+1}))^\gamma}
            \right)
        \frac{\rho_\Delta(x,t^{n+1}) dx}{\dx \rho_i^{n+1}}
        \right).
    \end{equation*}
    Since $\Xi^{-1}$ is a concave function,
    applying Jensen's inequality to the right-hand side of the above inequality, we obtain
    \begin{equation*}
        \frac{p(W_i^{n+1})}{(\rho_i^{n+1})^\gamma} \ge
        \int_{x_\imh}^{x_\iph}
        \frac{p(W_\Delta(x,t^{n+1}))}{(\rho_\Delta(x,t^{n+1}))^\gamma}
        \frac{\rho_\Delta(x,t^{n+1})}{\dx \rho_i^{n+1}}\:dx.
    \end{equation*}
    Since $W_\Delta(x,t^{n+1})$ is a juxtaposition of \ARSs defined by
    \eqref{eq:ARS}, because of the positivity of $\rho_i^{n+1}$ and
    $\rho_\Delta(x,t^{n+1})$ given by both $\rho^*_L>0$ and
    $\rho^*_R>0$, it results that $p(W_i^{n+1})>0$ as soon as
    \begin{equation*}
        p(W^*_L)>0
        \text{\qquad and \qquad}
        p(W^*_R)>0.
    \end{equation*}
    Finally, according to the definition \eqref{eq:entropy} of $s$,
    the pressures of the intermediate states satisfy
    $p(W^*_L) = (\rho_L^*)^\gamma \exp(-s^*)$
    and
    $p(W^*_R) = (\rho_R^*)^\gamma \exp(-s^*)$.
    Thus, the above positivity condition always holds,
    and $p(W_i^{n+1})$ is therefore positive.
\end{proof}

Equipped with these definitions
and conditions \eqref{eq:WB_consequence_on_hat_and_delta_w} --
\eqref{eq:WB_consequence_on_intermediate_states} --
\eqref{eq:condition_on_s_star} --
\eqref{eq:positivity_condition_on_rho},
we now have all the ingredients we need to determine $\delta W$
in \cref{sec:determination_of_delta_w},
and $S^q, S^E$ in \cref{sec:source_term_approximation}.

\subsection{\texorpdfstring{Determination of $\delta W$}{Determination of delta W}}
\label{sec:determination_of_delta_w}

The goal of this section is to use the definitions and conditions established in \cref{sec:definitions}
to derive an expression of the three yet unknown components of $\delta W$,
such that $\delta W = [W] / 2$ as soon as $(W_L, W_R)$ is an ISS,
where we have denoted by $[X] \coloneqq X_R - X_L$ the jump of some quantity $X$.
Moreover, for stability purposes,
we wish to recover the HLL solver in the absence of a gravitational source term.
Therefore, we also impose that,
if $\varphi_L = \varphi_R$, then $\delta W$ must vanish.

We start with the momentum component $\delta q$.
Note that, according
to~\eqref{eq:WB_consequence_on_intermediate_states},
we set $q_L^* = q_R^*$.
Therefore, $\delta q = 0$ is a suitable choice.

Now, we turn to the density component $\delta \rho$.
Recall that, if $(W_L, W_R)$ is an ISS, then $H(W_L, \varphi_L) = H(W_R, \varphi_R)$,
that is to say
\begin{equation}
    \label{eq:total_enthalpy_left_right_steady}
    \frac{E_L + p_L}{\rho_L} + \varphi_L
    =
    \frac{E_R + p_R}{\rho_R} + \varphi_R.
\end{equation}
Defining the enthalpy
\begin{equation}
    \label{eq:enthalpy}
    h = \frac{E + p}{\rho},
\end{equation}
and setting $h_L = h(W_L)$ and $h_R = h(W_R)$,
\eqref{eq:total_enthalpy_left_right_steady} rewrites as
\begin{equation}
    \label{eq:steady_relation_on_enthalpy}
    [h] = - [\varphi].
\end{equation}
In other words, if $[h] \neq 0$, then
\begin{equation}
    \label{eq:def_y}
    y \coloneqq 1 + \frac{[\varphi]}{[h]} = 0.
\end{equation}
The function $y \mapsto \delta \rho(y)$ must then fulfil the following properties:
\begin{align}
    \label{eq:deltarho_i} \tag{i}
    \delta \rho(0)                         & = [\rho]/2,  \\
    \label{eq:deltarho_ii} \tag{ii}
    \delta \rho(1)                         & = 0,         \\
    \tag{iii} \label{eq:deltarho_iii}
    |\delta \rho|                          & < \rho_\HLL, \\
    \tag{iv}  \label{eq:deltarho_iv}
    \lim_{y \to \pm \infty} \delta \rho(y) & = 0.
\end{align}
Property \eqref{eq:deltarho_i} corresponds to the ISS case,
and ensures that $\delta \rho$ satisfies
\eqref{eq:WB_consequence_on_hat_and_delta_w}.
Property \eqref{eq:deltarho_ii} corresponds to the no-gravity case,
since $y = 1$ as soon as $[\varphi] = 0$.
In this case, we wish $\delta \rho$ to vanish,
so that $\rho_L^* = \rho_R^* = \rho_\text{HLL}$.
Property \eqref{eq:deltarho_iii} is the required condition
\eqref{eq:positivity_condition_on_rho} for positivity preservation.
Property \eqref{eq:deltarho_iv} ensures that $\delta \rho$ is well-defined
and goes to $0$ when $[h]$ vanishes.
To determine the final expression of $\delta \rho$
satisfying properties \eqref{eq:deltarho_i} -- \eqref{eq:deltarho_iv},
we first introduce a relevant function $\psi$.

\begin{lemma}
    \label{lem:def_psi}
    Let $\psi$ be a smooth function given by
    \begin{equation}
        \label{eq:def_psi}
        \psi(y) = \cos \left(\frac \pi 2 y\right) \exp \left(-2 y^2 \right).
    \end{equation}
    Then $\psi$ satisfies the following properties:
    \begin{align}
        \label{eq:psi_i} \tag{$\psi$-i}
        \psi(0)                                 & =
        1,                                             \\
        \label{eq:psi_ii} \tag{$\psi$-ii}
        \psi(1)                                 & =
        0,                                             \\
        \label{eq:psi_iii} \tag{$\psi$-iii}
        |\psi|                                  & \leq
        1,                                             \\
        \label{eq:psi_iv} \tag{$\psi$-iv}
        \smash{\lim_{y \to \pm \infty}} \psi(y) & =
        0,                                             \\
        \label{eq:psi_v} \tag{$\psi$-v}
        \psi(1 + y)                             & =
        \mathcal{O}(y) \text{\qquad as \, $y \to 0$}.
    \end{align}
\end{lemma}

\begin{proof}
    Properties \eqref{eq:psi_i} and \eqref{eq:psi_ii}
    are straightforward to verify.
    For illustration, we provide in \cref{fig:psi}
    a graphical representation of $\psi$.
    Regarding property \eqref{eq:psi_iii}, we note that for all $y \in \mathbb{R}$,
    \begin{equation*}
        |\psi(y)|
        =
        \left| \cos \left(\frac \pi 2 y\right) \exp \left(-2 y^2 \right) \right|
        \leq \exp(-2 y^2) \leq 1.
    \end{equation*}
    Property \eqref{eq:psi_iv} holds since the cosine function is bounded
    and $y \mapsto e^{-2 y^2}$ goes to $0$ as~$y$ goes to $\pm \infty$.
    Finally, to prove property \eqref{eq:psi_v},
    we compute $\psi(1+y)$ for $y \in \mathbb{R}$.
    We obtain
    \begin{equation*}
        \psi(1+y)
        =
        \cos \left(\frac \pi 2 + \frac \pi 2 y \right) \exp \left(-2 (1+y)^2 \right)
        =
        - \sin \left(\frac \pi 2 y \right) \exp \left(-2 (1+y)^2 \right).
    \end{equation*}
    Therefore,
    \begin{equation*}
        \psi(1+y) \sim - \frac \pi 2 e^{-2} y \text{\qquad as \, $y \to 0$},
    \end{equation*}
    which proves \eqref{eq:psi_v} and concludes the proof.
\end{proof}

\begin{figure}[tb]
    \centering
    \begin{tikzpicture}[scale=0.7]
        \begin{scope}[every node/.style={scale=1.4}]
            \begin{axis}[
                    axis lines = left,
                    enlarge x limits={abs=5pt},
                    enlarge y limits={abs=5pt},
                    xlabel = {$y$},
                    xlabel style={at={(ticklabel* cs:1.01)},anchor=west},
                    ylabel = {\rotatebox{270}{$\psi$}},
                    ylabel style={at={(ticklabel* cs:1.01)},anchor=west},
                    xmin=-1.75, xmax=1.75,
                    ymin=-0.02, ymax=1.05,
                    domain=-1.75:1.75, samples=100,
                    hide obscured x ticks=false,
                    hide obscured y ticks=false,
                    scaled y ticks = false,
                    xtick={0,1},
                    ytick={0,1},
                    yticklabels={$0$, $1$},
                    grid=major,
                ]
                \addplot [smooth] {cos(90*x)*exp(-2*x^2)};
            \end{axis}
        \end{scope}
    \end{tikzpicture}
    \caption{%
        Graph of the function $\psi$ defined in \cref{lem:def_psi}.
        We observe that $\psi$ indeed satisfies
        properties \eqref{eq:psi_i} -- \eqref{eq:psi_iv}.
    }
    \label{fig:psi}
\end{figure}

We now have all the ingredients we need to construct $\delta \rho$ as a function of $y$.
Note that multiple expressions of such a function $y \mapsto \delta \rho(y)$ are possible.
We state our choice in the next Lemma.

\begin{lemma}
    \label{lem:positivity_density}
    Let $W_L \in \Omega$ and $W_R \in \Omega$ two given states,
    and let $W_\HLL$ be given by \eqref{eq:def_w_hll} with $\lambda$ sufficiently large.
    Further, let $y \mapsto \delta \rho(y)$ be a smooth function given by
    \begin{equation*}
        \label{eq:deltarho_definition}
        \delta \rho(y) =
        \frac{[\rho]} 2 \psi(y),
    \end{equation*}
    where $y$ is defined by \eqref{eq:def_y} and $\psi$ is given by \eqref{eq:def_psi}.
    Then, conditions \eqref{eq:deltarho_i} -- \eqref{eq:deltarho_iv}
    are satisfied by $\delta \rho(y)$.
\end{lemma}

\begin{proof}
    First, we note that the properties
    \eqref{eq:psi_i}, \eqref{eq:psi_ii} and \eqref{eq:psi_iv}
    directly map to properties \eqref{eq:deltarho_i},
    \eqref{eq:deltarho_ii} and \eqref{eq:deltarho_iv}.
    Second, to prove property \eqref{eq:deltarho_iii},
    we need to show that $|\delta \rho| < \rho_\HLL$.
    To that end, we analyze the extrema of the function~$\delta\rho(y)$.
    According to property \eqref{eq:psi_iii} of $\psi$,
    the global extremum of $\delta\rho(y)$ is $[\rho]/2$, reached for $y = 0$.
    Therefore, since $\rho_L > 0$ and $\rho_R > 0$,
    we immediately obtain
    \begin{equation*}
        |\delta \rho| \leq
        \frac{ |\rho_R - \rho_L| } 2 \leq
        \frac{\rho_L + \rho_R} 2.
    \end{equation*}
    Hence, according to the expression \eqref{eq:def_w_hll}
    of $\rho_\HLL$, there exists a large enough $\lambda$
    such that $|\delta \rho| < \rho_\HLL$,
    and property \eqref{eq:deltarho_iii} is satisfied.
    The proof is thus concluded.
\end{proof}

Finally, we turn to the energy component $\delta E$.
Analogously to the derivation of $\delta q$ and $\delta\rho$, we use the last remaining ISS condition
on $s$ to determine $\delta E$.
Namely, we wish to impose that $s(W_L^*) = s(W_R^*)$
as soon as $(W_L, W_R)$ is an ISS,
as prescribed in~\eqref{eq:WB_consequence_on_intermediate_states}.
Note that this implies, using the
definition \eqref{eq:entropy} of $s$, that
\begin{equation*}
    \label{eq:equality_of_intermediate_entropies}
    -\log\left( \frac{p(W_L^*)}{(\rho_L^*)^\gamma} \right)
    =
    -\log\left( \frac{p(W_R^*)}{(\rho_R^*)^\gamma} \right)
    \text{, \quad thus, \quad}
    \frac{p(W_L^*)}{(\rho_L^*)^\gamma}
    =
    \frac{p(W_R^*)}{(\rho_R^*)^\gamma}.
\end{equation*}
Therefore, using the expression \eqref{eq:pressure_law} of the pressure law, we obtain
\begin{equation*}
    \label{eq:Entropy_condition_deltaE_with_q_star}
    \frac{1}{(\rho_L^*)^\gamma} \left( E_L^* - \frac{(q_L^*)^2}{2 \rho_L^*} \right)
    =
    \frac{1}{(\rho_R^*)^\gamma} \left( E_R^* - \frac{(q_R^*)^2}{2 \rho_R^*} \right).
\end{equation*}
However, recall that, for a steady solution, the momentum is constant.
Hence, we replace $\smash{(q_L^*)^2}$ and $\smash{(q_R^*)^2}$ in the above relation
by a consistent average of $q^2$, denoted by $\smash{\widetilde{q^2}}$,
to be determined later on.
We then obtain
\begin{equation*}
    \label{eq:Entropy_condition_deltaE}
    \frac{1}{(\rho_L^*)^\gamma} \left( E_L^* - \frac{\widetilde{q^2}}{2 \rho_L^*} \right)
    =
    \frac{1}{(\rho_R^*)^\gamma} \left( E_R^* - \frac{\widetilde{q^2}}{2 \rho_R^*} \right),
\end{equation*}
Using $E_L^* = \widehat E - \delta E$ and $E_R^* = \widehat E + \delta E$ leads to
the following expression of $\delta E$
\begin{equation}
    \label{eq:expression_delta_E}
    \delta E =
    \frac 1 {(\rho_L^*)^\gamma + (\rho_R^*)^\gamma} \left(
    (\rho_R^*)^\gamma \left(\widehat E - \frac{\widetilde{q^2}}{2 \rho_L^*}\right)
    -
    (\rho_L^*)^\gamma \left(\widehat E - \frac{\widetilde{q^2}}{2 \rho_R^*}\right)
    \right).
\end{equation}
Note that this expression of $\delta E$ satisfies both properties
required for well-ba\-lan\-ced\-ness,
independently of the expression of $\smash{\widetilde{q^2}}$.
Hence, $\smash{\widetilde{q^2}}$ is a degree of freedom in the numerical scheme.
It turns out that
one can determine $\smash{\widetilde{q^2}}$
such that the scheme is entropy-satisfying,
by leveraging condition \eqref{eq:condition_on_s_star}
on the intermediate entropy~$s^\ast$.
Indeed, recall that we have defined $s^* = s(W_L^*) = s(W_R^*)$,
and that $s^*$ is given by $s^* = (\rho s)_\HLL / \rho_\HLL$.
By definition of a constant-entropy intermediate state, we obtain
\begin{equation*}
    -\log \left( \frac{p(W_L^*)}{(\rho_L^*)^\gamma} \right)
    =
    s^*
    =
    \frac{(\rho s)_\HLL}{\rho_\HLL}.
\end{equation*}
Therefore under application of the exponential function, we have
\begin{equation*}
    \frac 1 {(\rho_L^*)^\gamma}
    \left( \widehat E - \delta E - \frac{\widetilde{q^2}}{2 \rho_L^*} \right)
    =
    \frac 1 {\gamma - 1}
    \exp \left( -\frac{(\rho s)_\HLL}{\rho_\HLL} \right).
\end{equation*}
Substituting $\delta E$ by  \eqref{eq:expression_delta_E} in the above relation leads,
after straightforward computations, to
\begin{equation*}
    \label{eq:expression_tilde_q2}
    \widetilde{q^2} =
    \frac {2 \rho_L^* \rho_R^*} {\rho_L^* + \rho_R^*}
    \left(
    2 \widehat E
    -
    \exp \left( -\frac{(\rho s)_\HLL}{\rho_\HLL} \right)
    \frac{(\rho_L^*)^\gamma + (\rho_R^*)^\gamma}{\gamma - 1}
    \right).
\end{equation*}
Note that this expression is indeed consistent with a squared momentum,
since the expression in brackets is consistent
with a kinetic energy.


\subsection{\texorpdfstring{Determination of $S^q$ and $S^E$}{Determination of Sq and SE}}
\label{sec:source_term_approximation}

Equipped with the expression of $\delta W$,
we now turn to the determination of the last two unknowns, $S^q$ and $S^E$.
To that end, we leverage the expressions \eqref{eq:expression_w_hat} of $\widehat W$,
and recall that, if $(W_L, W_R)$ is an ISS, then $\widehat W = \xoverline W$,
where we have defined $\xoverline X \coloneqq (X_L + X_R) / 2$
the arithmetic mean of two quantities $X_L$ and $X_R$.
Therefore, plugging the expression of $W_\HLL$ in \eqref{eq:expression_w_hat},
we obtain that, if $(W_L, W_R)$ is an ISS, then $S^q$ and $S^E$ must satisfy
\begin{subequations}
    \label{eq:conditions_on_Sq_and_SE}
    \begin{align}
        \label{eq:conditions_on_Sq}
        S^q \dx & =
        \left(\frac{q_R^2}{\rho_R} + p_R\right)
        -
        \left(\frac{q_L^2}{\rho_L} + p_L\right), \\
        \label{eq:conditions_on_SE}
        S^E \dx & =
        \left(\frac{q_R}{\rho_R} (E_R + p_R) \right)
        -
        \left(\frac{q_L}{\rho_L} (E_L + p_L)\right).
    \end{align}
\end{subequations}
Moreover, recall that $S^q$ and $S^E$ should be consistent approximations
of the source term averages.
Therefore, the goal of this section is to find such consistent averages
satisfying \eqref{eq:conditions_on_Sq_and_SE}.

\begin{remark}
    \label{rem:smoothness_of_potential}
    Since $\varphi$ is assumed to be a smooth function,
    for each interface, the following identities hold:
    \begin{equation*}
        \varphi_R = \varphi_L + \mathcal{O}(\dx)
        \text{, \quad i.e., \quad}
        [\varphi] = \mathcal{O}(\dx).
    \end{equation*}
\end{remark}

We start with the energy source term $S^E$.
Its expression and related properties are given in the following lemma.

\begin{lemma}
    Let two given states
    $W_L \in \Omega$ and $W_R \in \Omega$.
    Further, let
    \begin{equation}
        \label{eq:expression_of_SE}
        S^E = - \frac{q_L + q_R} 2 \, \frac{\varphi_R - \varphi_L} \dx.
    \end{equation}
    Then $S^E$ is consistent with $-q \pd_x \varphi$
    and satisfies the well-balanced condition \eqref{eq:conditions_on_SE}.
\end{lemma}

\begin{proof}
    To prove consistency with the continuous source term,
    we take $q_L = q(x)$ and $q_R = q(x + \dx)$,
    with $\varphi_L = \varphi(x)$ and $\varphi_R = \varphi(x + \dx)$.
    Using \cref{rem:smoothness_of_potential},
    a simple Taylor expansion of $S^E$ then shows that
    $S^E = - q \pd_x \varphi + \mathcal{O}(\dx)$,
    which proves the consistency of $S^E$.

    For the well-balancedness,
    recall that the enthalpy is given
    in \eqref{eq:enthalpy} and satisfies $E + p = \rho h$.
    So, if $(W_L, W_R)$ is an ISS with $q_L = q_R = \xoverline q$,
    then $S^E \dx$ must satisfy
    $S^E \dx = \xoverline q (h_R - h_L)$ after \eqref{eq:conditions_on_SE}.
    However, from \eqref{eq:steady_relation_on_enthalpy} we have
    that $[h] = - [\varphi]$ if $(W_L, W_R)$ is an ISS.
    Therefore, $S^E$ given by \eqref{eq:expression_of_SE}
    satisfies \eqref{eq:conditions_on_SE} by construction,
    which concludes the proof.
\end{proof}

Next, we consider the momentum source term $S^q$.
Its expression and derivation are more involved and are summarized in the following lemma.

\begin{lemma}
    Let two given states
    $W_L \in \Omega$ and $W_R \in \Omega$.
    Further, let
    \begin{equation}
        \label{eq:expression_of_Sq}
        S^q
        =
        - \frac{2 \rho_L \rho_R}{\rho_L + \rho_R} \, \frac{\varphi_R - \varphi_L}{\dx}
        + \frac \varepsilon \dx \, \psi \left(
        1 + \left( \frac{\varphi_R - \varphi_L}{h_R - h_L} \right)^3
        \right),
    \end{equation}
    where $\psi$ is given by \eqref{eq:def_psi} and $\varepsilon$ is given by
    \begin{equation}
        \label{eq:expression_of_epsilon}
        \varepsilon = \xoverline{e^{-s}} \left(
        \left( \rho_R^\gamma - \rho_L^\gamma \right)
        -
        \frac{2 \rho_L \rho_R}{\rho_L + \rho_R} \,
        \frac{\gamma}{\gamma - 1} \,
        \left( \rho_R^{\gamma - 1} - \rho_L^{\gamma - 1} \right)
        \right),
    \end{equation}
    with $\xoverline{e^{-s}} = (e^{-s_L} + e^{-s_R}) / 2$.
    Then $S^q$ is consistent with $-\rho \pd_x \varphi$
    and satisfies the well-balanced condition \eqref{eq:conditions_on_Sq}.
\end{lemma}

\begin{proof}
    We first prove that $S^q$ is consistent with the continuous source term.
    In this context, we emphasize that $(W_L, W_R)$ is not necessarily an ISS,
    but rather a general pair of states.
    Performing a Taylor expansion of the first term of~$S^q$,
    and recalling \cref{rem:smoothness_of_potential}, we immediately see that
    \begin{equation*}
        - \frac{2 \rho_L \rho_R}{\rho_L + \rho_R} \, \frac{\varphi_R - \varphi_L}{\dx}
        \equals_{\dx \to 0} - \rho \pd_x \varphi + \mathcal{O}(\dx).
    \end{equation*}
    For the second term of $S^q$, we distinguish two cases:
    either the solution is smooth and $W_R = W_L + \mathcal{O}(\dx)$,
    or the solution is not smooth (e.g. in a shock wave) and $W_R = W_L + \mathcal{O}(1)$.
    It turns out that the multiplication by $\psi$ in \eqref{eq:expression_of_Sq}
    is crucial to overcome a consistency defect in the non-smooth case.
    In the first case, performing a Taylor expansion of $\varepsilon$
    given by~\eqref{eq:expression_of_epsilon}
    and noting that $[\varphi]/[h] = \mathcal{O}(1)$
    since both~$[\varphi]$ and $[h]$ are $\mathcal{O}(\dx)$,
    we obtain that
    \begin{equation*}
        W_R = W_L + \mathcal{O}(\dx)
        \implies
        \varepsilon = \mathcal{O}(\dx^3)
        \text{\quad and \quad}
        \psi\left(1 + \frac{[\varphi]^3}{[h]^3}\right) = \mathcal{O}(1).
    \end{equation*}
    In the second case, using property \eqref{eq:psi_v} of $\psi$
    and remarking that $[\varphi] = \mathcal{O}(\dx)$ and $[h] = \mathcal{O}(1)$, we get
    \begin{equation*}
        W_R = W_L + \mathcal{O}(1)
        \implies
        \varepsilon = \mathcal{O}(1)
        \text{\quad and \quad}
        \psi\left(1 + \frac{[\varphi]^3}{[h]^3}\right)
        =
        \mathcal{O}\left(\frac{[\varphi]^3}{[h]^3}\right)
        =
        \mathcal{O}(\dx^3).
    \end{equation*}
    As a consequence, in both cases, the second term of $S^q$ is consistent with $0$,
    with an error of $\mathcal{O}(\dx^2)$.
    Therefore, we have obtained that $S^q$ is
    indeed consistent with the continuous source term,
    i.e., that $S^q = -\rho \pd_x \varphi + \mathcal{O}(\dx)$.

    We now turn to the well-balanced property;
    we have to prove that $S^q$, given by~\eqref{eq:expression_of_Sq},
    satisfies \eqref{eq:conditions_on_Sq} as soon as $(W_L, W_R)$ is an ISS.
    In this case, we recall that $[\varphi] = -[h]$.
    Therefore, using property \eqref{eq:psi_i} of $\psi$,
    we obtain $\psi(1 + ([\varphi]/[h])^3) = 1$.
    Hence, we have to prove that
    \begin{equation*}
        - \frac{2 \rho_L \rho_R}{\rho_L + \rho_R} \, (\varphi_R - \varphi_L)
        + \varepsilon
        =
        \left(\frac{q_R^2}{\rho_R} + p_R\right)
        -
        \left(\frac{q_L^2}{\rho_L} + p_L\right)
    \end{equation*}
    as soon as $(W_L, W_R)$ is an ISS
    and for $\varepsilon$ given by \eqref{eq:expression_of_epsilon}.
    To that end, we first address $\varepsilon$.
    Note that, since we are considering an ISS, we have $s_L = s_R$.
    Therefore,
    \begin{equation*}
        \xoverline{e^{-s}}
        =
        \frac 1 2 \left(
        e^{-s_L} + e^{-s_R}
        \right)
        =
        \frac 1 2 \left(
        \frac{p_L} {\rho_L^\gamma} + \frac{p_R} {\rho_R^\gamma}
        \right)\!
        =
        \frac{p_L} {\rho_L^\gamma}
        =
        \frac{p_R} {\rho_R^\gamma}.
    \end{equation*}
    Using these identities to distribute $\xoverline{e^{-s}}$, $\varepsilon$ rewrites
    \begin{equation*}
        \varepsilon =
        \left( p_R - p_L \right)
        -
        \frac{2 \rho_L \rho_R}{\rho_L + \rho_R} \,
        \frac{\gamma}{\gamma - 1} \,
        \left( \frac{p_R}{\rho_R} - \frac{p_L}{\rho_L} \right),
    \end{equation*}
    and, from the expression \eqref{eq:expression_of_Sq} of $S^q$, we obtain
    \begin{equation}
        \label{eq:S_q_intermediate_expression}
        S^q \dx
        =
        \left( p_R - p_L \right)
        -
        \frac{2 \rho_L \rho_R}{\rho_L + \rho_R} \,
        \left(
        (\varphi_R - \varphi_L) +
        \frac{\gamma}{\gamma - 1} \,
        \left( \frac{p_R}{\rho_R} - \frac{p_L}{\rho_L} \right)
        \right).
    \end{equation}
    Now, using the steady relation $[\varphi] = -[h]$
    together with the definition \eqref{eq:enthalpy} of $h$,
    we get, using \eqref{eq:enthalpy_wrt_rho_q_p},
    \begin{equation*}
        \varphi_R - \varphi_L
        =
        - \left(\frac{E_R + p_R}{\rho_R} - \frac{E_L + p_L}{\rho_L}\right)
        =
        - \frac{\gamma}{\gamma - 1} \left(
        \frac{p_R}{\rho_R} - \frac{p_L}{\rho_L}
        \right)
        - \left(
        \frac 1 2 \frac{q_R^2}{\rho_R^2} - \frac 1 2 \frac{q_L^2}{\rho_L^2}
        \right).
    \end{equation*}
    Plugging this relation into \eqref{eq:S_q_intermediate_expression},
    the term in brackets reduces to
    \begin{equation*}
        (\varphi_R - \varphi_L) +
        \frac{\gamma}{\gamma - 1} \,
        \left( \frac{p_R}{\rho_R} - \frac{p_L}{\rho_L} \right)
        =
        - \left(
        \frac 1 2 \frac{q_R^2}{\rho_R^2} - \frac 1 2 \frac{q_L^2}{\rho_L^2}
        \right).
    \end{equation*}
    Now, recall that one of the ISS relations is $q_L = q_R = \xoverline q$.
    Consequently, we obtain
    \begin{equation*}
        \frac{2 \rho_L \rho_R}{\rho_L + \rho_R}
        \left( \frac 1 2 \frac{q_R^2}{\rho_R} - \frac 1 2 \frac{q_L^2}{\rho_L} \right)
        =
        \xoverline q^2 \frac{2 \rho_L \rho_R}{\rho_L + \rho_R}
        \left( \frac{1}{\rho_R^2} - \frac{1}{\rho_L^2} \right)
        =
        \xoverline q^2
        \left( \frac{1}{\rho_R} - \frac{1}{\rho_L} \right)
        =
        \frac{q_R^2}{\rho_R} - \frac{q_L^2}{\rho_L}.
    \end{equation*}
    Substituting into \eqref{eq:S_q_intermediate_expression} leads to
    \begin{equation*}
        S^q \dx
        =
        \left( p_R - p_L \right)
        +
        \left( \frac{q_R^2}{\rho_R} - \frac{q_L^2}{\rho_L} \right),
    \end{equation*}
    which is nothing but the well-balancedness condition \eqref{eq:conditions_on_SE}.
    The proof is thus concluded.
\end{proof}

\subsection{Summary of the numerical scheme and main properties}

We finally summarize the approximate Riemann solver
and state the main properties of the scheme deriving from this \ARS.
The \ARS \eqref{eq:ARS} is given by
\begin{equation}
    \label{eq:ARS_recap}
    \widetilde{W}\left( \frac x t; W_L, W_R \right) =
    \begin{dcases}
        W_L                           & \text{if } x < -\lambda t,     \\
        W_L^* = \widehat W - \delta W & \text{if } -\lambda t < x < 0, \\
        W_R^* = \widehat W + \delta W & \text{if } 0 < x < \lambda t,  \\
        W_R                           & \text{if } x > \lambda t.      \\
    \end{dcases}
\end{equation}
This \ARS is based on the HLL solver,
whose intermediate state is given by \eqref{eq:def_w_hll},
with the wave speed $\lambda$ given by \eqref{eq:def_lambda}.
Recall also the pressure law \eqref{eq:pressure_law} defining $p$,
the expression \eqref{eq:entropy} of the entropy $s$,
and the expression \eqref{eq:enthalpy_wrt_rho_q_p} of the enthalpy $h$.
Equipped with these definitions,
\eqref{eq:expression_w_hat} gives the following expression of $\smash{\widehat W}$:
\begin{equation}
    \label{eq:expression_w_hat_recap}
    \left\{
    \begin{aligned}
        \widehat \rho & = \rho_\HLL,   \\
        \widehat q    & = q_\HLL
        + \frac {S^q \dx} {2 \lambda}, \\
        \widehat E    & = E_\HLL
        + \frac {S^E \dx} {2 \lambda}, \\
    \end{aligned}
    \right.
\end{equation}
where the approximate source terms $S^q$ and $S^E$
have been derived in \cref{sec:source_term_approximation}
and are given by
\begin{equation}
    \begin{aligned}
        \label{eq:expressions_of_Sq_and_SE_recap}
        S^E & = - \frac{q_L + q_R} 2 \, \frac{\varphi_R - \varphi_L} \dx, \\
        S^q
            & =
        - \frac{2 \rho_L \rho_R}{\rho_L + \rho_R} \, \frac{\varphi_R - \varphi_L}{\dx}
        + \frac \varepsilon \dx \, \psi \left(
        1 + \left( \frac{\varphi_R - \varphi_L}{h_R - h_L} \right)^3
        \right),
    \end{aligned}
\end{equation}
where $\psi(y) = \cos(\frac \pi 2 y) e^{-2y^2}$,
and where $\varepsilon$ satisfies
\begin{equation*}
    \label{eq:expression_of_epsilon_recap}
    \varepsilon = \frac 1 2 \left( e^{-s_L} + e^{-s_R} \right) \left(
    \left( \rho_R^\gamma - \rho_L^\gamma \right)
    -
    \frac{2 \rho_L \rho_R}{\rho_L + \rho_R} \,
    \frac{\gamma}{\gamma - 1} \,
    \left( \rho_R^{\gamma - 1} - \rho_L^{\gamma - 1} \right)
    \right).
\end{equation*}
The components of $\delta W$ have been derived in \cref{sec:determination_of_delta_w},
and we have set
\begin{equation}
    \label{eq:expression_delta_w_recap}
    \left\{
    \begin{aligned}
        \delta \rho &
        = \frac {\rho_R - \rho_L} 2 \,
        \psi \! \left(1 + \frac{\varphi_R - \varphi_L}{h_R - h_L} \right), \\
        \delta q    & = 0 \vphantom{\dfrac 1 2},                           \\
        \delta E    & =
        \frac 1 {(\rho_L^*)^\gamma + (\rho_R^*)^\gamma} \left(
        (\rho_R^*)^\gamma \left(\widehat E - \frac{\widetilde{q^2}}{2 \rho_L^*}\right)
        -
        (\rho_L^*)^\gamma \left(\widehat E - \frac{\widetilde{q^2}}{2 \rho_R^*}\right)
        \right),                                                           \\
    \end{aligned}
    \right.
\end{equation}
where $\widetilde{q^2}$ satisfies
\begin{equation*}
    \label{eq:expression_tilde_q2_recap}
    \widetilde{q^2} =
    \frac {2 \rho_L^* \rho_R^*} {\rho_L^* + \rho_R^*}
    \left(
    2 \widehat E
    -
    \exp \left( -\frac{(\rho s)_\HLL}{\rho_\HLL} \right)
    \frac{(\rho_L^*)^\gamma + (\rho_R^*)^\gamma}{\gamma - 1}
    \right).
\end{equation*}

Equipped with these expressions, we are ready to state our main result,
summarizing all properties of the resulting scheme.

\begin{theorem}
    \label{theo:summary}
    Let the time step $\dt$ be given by \eqref{eq:CFL_condition}
    and assume that the initial data satisfies
    $W_i^0 \in \Omega$ for all $i \in \mathbb{Z}$.
    Then, the numerical scheme \eqref{eq:scheme_from_ARS}
    with the approximate Riemann solver \eqref{eq:ARS_recap},
    where $\smash{\widehat{W}}$ is given by \eqref{eq:expression_w_hat_recap}
    and $\delta W$ is given by \eqref{eq:expression_delta_w_recap},
    satisfies the following properties:
    \begin{enumerate}
        \item consistency with the Euler system \eqref{eq:EulerG};
        \item positivity preservation: for all $n \geq 0$,
              \begin{equation*}
                  \forall i \in \mathbb{Z}, \; W_i^{n} \in \Omega
                  \quad \implies \quad
                  \forall i \in \mathbb{Z}, \; W_i^{n+1} \in \Omega;
              \end{equation*}
        \item entropy stability:
              for all smooth functions $\eta$ satisfying \eqref{eq:entropy_cond},
              for all $i \in \mathbb{Z}$, for all $n \geq 0$,
              \begin{equation*}
                  \rho_i^{n+1}\eta(s_i^{n+1}) \leq \rho_i^n \eta(s_i^n)
                  - \frac{\dt}{\dx} \left((\rho \eta(s) u)_{i+1/2}^n -
                  (\rho \eta(s) u)_{i-1/2}^n\right);
              \end{equation*}
        \item well-balancedness:
              \begin{equation*}
                  \forall i \in \mathbb{Z}, \; (W_i^n, W_{i+1}^n) \text{ is an ISS}
                  \quad \implies \quad
                  \forall i \in \mathbb{Z}, \; W_i^{n+1} = W_i^n.
              \end{equation*}
    \end{enumerate}
\end{theorem}

\begin{proof}
    We prove the four properties in order,
    using the results derived in the previous sections.
    \begin{enumerate}
        \item According to \cite{HarLaxLee1983},
              the scheme is consistent as soon as the \ARS satisfies
              the integral consistency relation \eqref{eq:integral_consistency},
              which holds by construction,
              as detailed in \cref{sec:Riemann_solver_consistency}.
              Thus, the \ARS is consistent,
              and therefore the numerical scheme \eqref{eq:scheme_from_ARS}
              is also consistent.
        \item We have to prove the positivity
              of the density and the pressure.
              A sufficient condition for the positivity of the intermediate density
              is given by \eqref{eq:positivity_condition_on_rho}.
              It is satisfied due to \cref{lem:positivity_density},
              which yields $\rho_L^\ast > 0$ and $\rho_R^\ast > 0$ in the \ARS,
              and thus $\rho_i^{n+1} > 0$ in the numerical scheme.
              As a consequence, \cref{lem:positivity_pressure} can be applied,
              yielding the positivity of the pressure $p_i^{n+1} > 0$.
        \item After \cite{HarLaxLee1983}, the entropy stability holds
              as soon as \eqref{eq:entropy_integral_consistency} is satisfied.
              However, following \cref{th:entrop},
              this relation is verified if the intermediate entropy
              of the \ARS is given by \eqref{eq:condition_on_s_star},
              which is the case for the \ARS \eqref{eq:ARS_recap}.
              This proves that the scheme is entropy-stable.
        \item A sufficient condition for well-balancedness
              is given in \cref{cor:WB_consequence},
              from which follow the sufficient conditions
              \eqref{eq:WB_consequence_on_hat_and_delta_w}
              to be satisfied by $\smash{\widehat{W}}$ and $\delta W$
              in the intermediate states of the \ARS.
              As shown in \cref{sec:determination_of_delta_w,sec:source_term_approximation},
              these conditions are satisfied by construction,
              and consequently the numerical scheme is fully well-balanced.
    \end{enumerate}
    Thus, all four properties have been proven, which concludes the proof.
\end{proof}

\section{A well-balanced higher order extension}
\label{sec:HigherOrder}

Equipped with the first-order scheme described
in \cref{theo:summary}, we now turn to the construction
of a high-order well-balanced extension.
Since the first-order scheme does not require solving
any non-linear system of equations,
we wish to preserve this property for the high-order scheme.
Other strategies have been developed over the years
to produce high-order well-balanced schemes,
but most of them require solving non-linear systems, see e.g. the non-exhaustive
list~\cite{CasGalLopPar2008,Xin2014,BriXin2020,GomCasParRus2021}.
To avoid an ill-posed problem, i.e. multiple possible solutions
associated to solving non-linear systems,
we follow the procedure described in detail in~\cite{BerBulFouMbaMic2022}
and used e.g. in ~\cite{MicBerClaFou2021}
for the shallow water equations.
Note that by adopting a strategy that circumvents non-linear iterative solvers, the associated computational overhead is avoided.

Here, we briefly outline the chosen procedure
to ensure high-order accuracy in space.
High-order accuracy in time will be obtained
by using traditional Runge-Kutta-type solvers;
it is described at the beginning of \cref{sec:NumRes}.
Blending high-order accuracy with the well-balanced property
is based on the observation that,
if the solution is steady,
the well-balanced scheme is exact,
and as such it will be more accurate than any high-order scheme.
As a consequence, this procedure consists in computing
the error between the solution and a steady solution,
and using this error to detect a steady solution
and correct the high-order scheme.
First, we introduce the uncorrected high-order scheme, see \cite{DioClaLou2012} for further details.
For a given $d \geq 1$, the scheme of order $d+1$ in space is
based on the polynomial reconstruction
\begin{equation}
    \label{eq:polynomial_reconstruction}
    \widecheck{W}_i^n(x) = W_i^n + \Pi_i^n(x - x_i),
\end{equation}
where $\Pi_i^n$ is a polynomial of degree $d$,
defined such that the reconstruction satisfies the following two relations
\begin{equation*}
    \label{eq:polynomial_reconstruction_properties}
    \widecheck{W}_i^n(x) = W(x, t^n) + \mathcal{O}(\dx^{d+1})
    \text{\quad and \quad}
    \frac 1 \dx \int_{x_\imh}^{x_\imh} \widecheck{W}_i^n(x) \, dx = W_i^n.
\end{equation*}
At this point we want to remark that it is at times more convenient to reconstruct in primitive variables $\rho$, $u$, $p$ than in conservative variables $W$.
Equipped with the reconstruction
\eqref{eq:polynomial_reconstruction},
interface reconstructions on the cell $C_i$ are given by
\begin{equation*}
    \label{eq:high_order_interface_reconstructions}
    \widecheck{W}_{\iph, -}^n = \widecheck{W}_i^n(x_\iph)
    \text{\quad and \quad}
    \widecheck{W}_{\iph, +}^n = \widecheck{W}_{i+1}^n(x_\iph),
\end{equation*}
while the high-order source term satisfies
\begin{equation}
    \label{eq:high_order_source_term}
    \widecheck{S}_i^n =
    \frac 1 \dx \int_{x_\imh}^{x_\imh} S(W(x, t^n)) \, dx
    + \mathcal{O}(\dx^{d+1}).
\end{equation}
The integral in \eqref{eq:high_order_source_term} is computed with
a quadrature rule of order $d$.
Then, using the interface values $\smash{\widecheck W_{i\pm1/2,\pm}^n}$
and the high-order source terms $\smash{\widecheck S_i}$, the scheme of order $d$ in space reads
\begin{equation}
    \label{eq:high_order_scheme}
    W_i^{n+1} = W_i^n
    -
    \frac{\dt}{\dx} \left(
    \mathcal{F}\big(\widecheck{W}_{\iph, -}^n, \widecheck{W}_{\iph, +}^n\big)
    -
    \mathcal{F}\big(\widecheck{W}_{\imh, -}^n, \widecheck{W}_{\imh, +}^n\big)
    \right)
    +
    \dt \widecheck{S}_i^n.
\end{equation}
As pointed out in \cref{sec:Godunov_type_scheme},
the numerical flux $\mathcal{F}$ defined in \eqref{eq:num_flux}
also implicitely depends on the gravity potentials at the
left and right of the interface.
In the context of \cref{sec:WBscheme},
these corresponded to cell averages.
In the current high-order context,
these left and right gravity potentials
are nothing but polynomial reconstructions of $\varphi$,
in cells $i$ and $i+1$, evaluated at \smash{$x_\iphinline$}.
Once again, we do not explicitly write them as inputs
of $\mathcal{F}$, to avoid cluttering notation.

Note that the scheme \eqref{eq:high_order_scheme}
is high-order accurate in space,
but not well-balanced.
Therefore, we employ the indicator from \cite{BerBulFouMbaMic2022}, to detect whether the current solution approximates a steady state.
It relies on the verification of the ISS property \cref{def:ISS} on all cells, thus the quantity
\begin{equation}
    \label{eq:sigma}
    \sigma_\iph^n
    =
    \left\|
    \begin{pmatrix}
        q_{i+1}^n - q_{i}^n       \\
        H(W_{i+1}^n, \varphi_{i+1}^n)
        -
        H(W_{i}^n, \varphi_{i}^n) \\
        s(W_{i+1}^n) - s(W_{i}^n) \\
    \end{pmatrix}
    \right\|.
\end{equation}
vanishes if, and only if,
\smash{$(W_i^n, W_{i+1}^n)$} is an ISS.
Thus \smash{$\sigma_\iphinline^n$} is well-suited as a steady solution detector.
To allow for a self adaptive correction, the high-order reconstructed interface values $\smash{\widecheck W_{i\pm1/2,\pm}^n}$ and source term discretization~$\smash{\widecheck S_i^n}$ are modified according to a parameter $\theta_\iphinline \in [0,1]$ yielding the cell average values in case of a steady solution, thus locally reverting to the first-order well-balanced scheme.
It is given by
\begin{equation}
    \label{eq:interface_recon_convex}
    \widebreve{W}_{\iph, \pm}^n =\left(1-  \theta_\iph^n\right)
    W_i^n + \theta_\iph^n\widecheck W_{i+1/2,\pm}
\end{equation}
Therein, the indicator is based on the detector \eqref{eq:sigma} and is defined by
\begin{equation}
    \label{eq:def_theta}
    \theta_\iph^n
    =
    \frac
    {\sigma_\iph^n}
    {\sigma_\iph^n + \left(\dfrac{\dx}{C_\iph^n}\right)^{d+1}}.
\end{equation}
Note that $\theta_\iphinline^n = 0$
as soon as $(W_i^n, W_{i+1}^n)$ is an ISS,
and otherwise yields
$\theta_\iphinline^n = 1 - \mathcal{O}(\smash{\dx^{d+1}}).$
Thus, the modified reconstruction given in \eqref{eq:interface_recon_convex}
satisfies
\begin{equation*}
    \widebreve{W}_{\iph, -}^n =
    \begin{cases}
        W_i^n
         & \text{ if } (W_i^n, W_{i+1}^n) \text{ is an ISS}, \\
        \widecheck{W}_i^n(x_\iph) + \mathcal{O}(\dx^{d+1})
         & \text{ otherwise},
    \end{cases}
\end{equation*}
which indeed shows that the scheme will be well-balanced.
In \eqref{eq:def_theta}, the quantity $C_\iph^n$ is given by
$C_\iph^0 = 0$ in the first time step $n=0$ and
\begin{equation}
    \label{eq:def_C_iph}
    C_\iph^n = C_\theta \,
    \frac 1 2
    \left(
    \frac { \|W_{i+1}^n - W_{i+1}^{n-1}\| } \dt
    +
    \frac { \|W_{i}^n - W_{i}^{n-1}\| } \dt
    \right)
    \text{ \quad if }n \geq 1.
\end{equation}
The reader is referred to \cite{BerBulFouMbaMic2022}
for a discussion regarding the expression of \smash{$C_\iphinline^n$},
and to the following section
for the choice of the parameter $C_\theta$ in the numerical experiments.

Similarly to \eqref{eq:interface_recon_convex}, the modified
high-order well-balanced source term is defined as
\begin{equation}
    \label{eq:high_order_WB_source_term}
    \widebreve{S}_i^n =
    \left( 1 - \frac {\theta_\imh^n + \theta_\iph^n} 2 \right) \widecheck{S}_i^n
    +
    \frac{\theta_\iph^n\mathcal{S}_\iph^n + \theta_\imh^n \mathcal{S}_\imh^n}{2}.
\end{equation}

Summarizing, using \eqref{eq:interface_recon_convex} and \eqref{eq:high_order_WB_source_term}, the scheme that is both well-balanced and high-order accurate in space reads
\begin{equation}
    \label{eq:high_order_WB_scheme}
    W_i^{n+1} = W_i^n
    -
    \frac{\dt}{\dx} \left(
    \mathcal{F}\big(\widebreve{W}_{\iph, -}^n, \widebreve{W}_{\iph, +}^n\big)
    -
    \mathcal{F}\big(\widebreve{W}_{\imh, -}^n, \widebreve{W}_{\imh, +}^n\big)
    \right)
    +
    \dt \widebreve{S}_i^n.
\end{equation}
To prevent spurious oscillations around discontinuities due to the high order reconstruction, the scheme \eqref{eq:high_order_WB_scheme} is endowed with a slope limiter.
Since this is merely modifying the reconstruction operator $\Pi_i^n$ in \eqref{eq:polynomial_reconstruction} the well-balanced property is not lost.

\section{Numerical results}
\label{sec:NumRes}

In this section, we perform numerical test cases to
validate the theoretical properties of the
first, second and third order well-balanced schemes,
denoted respectively by
$\mathbb{P}_0$, $\mathbb{P}_1^\WB$ and $\mathbb{P}_2^\WB$.
The $\mathbb{P}_1^\WB$ scheme
uses the \textsf{minmod} limiter, see for instance \cite{Lee1979},
and the SSPRK2 time integrator,
while the $\mathbb{P}_2^\WB$ scheme
makes uses the third order TVD reconstruction from \cite{SchSeiTor2015}
and the SSPRK3 time integrator.
Both Runge-Kutta schemes can be found
for instance in \cite{GotShuTad2001}.

For all test cases, we consider the ideal gas
EOS \eqref{eq:pressure_law} with $\gamma = 1.4$.
Further, we will use the gravitational potentials $\varphi_1(x) = (x-0.5)^2 / 2$ and $\varphi_2(x) = \sin(x)$.
Unless otherwise mentioned,
the parameters $\Lambda$ in \eqref{eq:def_lambda}
and $C_\theta$ in \eqref{eq:def_C_iph} are both set to $1$.

To check all properties of the scheme,
we validate the order of accuracy in \cref{sec:numerical_accuracy},
the well-balanced property in
\cref{sec:numerical_WB,sec:numerical_WB_pert}
on unperturbed and perturbed steady solutions respectively,
and the consistency with the Euler equations in
\cref{sec:numerical_RP}
by considering Riemann problems.
Lastly, a two-dimensional extension
is provided in \cref{sec:2D}.

\subsection{Accuracy of the numerical schemes}
\label{sec:numerical_accuracy}

\begin{figure}[tb]
    \centering
    \includegraphics[width=0.8\textwidth]{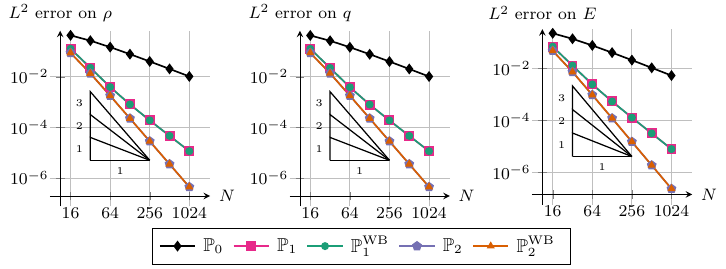}
    \caption{%
        Analytical solution \eqref{eq:exact_sol} with $\varphi(x) = x$:
        $L^2$ errors and EOC for the first, second and third order schemes,
        with or without the high-order well-balanced correction.%
    }
    \label{fig:EOC_phi1}
\end{figure}

To verify the experimental order of convergence (EOC) of the numerical schemes,
we study an exact solution for the Euler equations with gravity \eqref{eq:EulerG}
taken from~\cite{KliPupSem2019},
which is a variation of the test introduced in \cite{XinShu2012}.
The analytical solution is given by
\begin{equation}
    \label{eq:exact_sol}
    \begin{dcases}
        \rho(x,t)
        = 1 + \frac{1}{5} \sin(k \pi (x - u_0 t)), \\
        u(x,t)
        = u_0 \vphantom{\dfrac 1 5},               \\
        p(x,t)
        = \frac{9}{2} - ( x - u_0 t) + \frac{1}{k\pi} \cos(k\pi (x - u_0 t)),
    \end{dcases}
\end{equation}
where $u_0$ is a constant velocity, taken as $u_0 = 1$ in the experiment.
We take $k = 5$, which yields a highly oscillatory solution.
This test case is designed for a linear potential $\varphi(x) = x$.
The simulation is carried out on the computational domain $\Omega = [0,2]$
up to a final time of $t_f = 0.1$,
on seven grids with $N = 16 \cdot 2^{j}$ cells, where $j \in \{ 0, \dots, 6 \}$.
We prescribe exact boundary conditions.
In \cref{fig:EOC_phi1}, the $L^2$ errors are depicted
with respect to the number of cells.
We recover the expected order of convergence for all schemes.
Namely, we observe that the high-order well-balanced correction
does not negatively affect the EOC of the respective schemes.

\subsection{Well-balanced test cases}
\label{sec:numerical_WB}

To numerically verify the well-balanced property of the schemes,
we compute the numerical solutions of a hydrostatic and a moving equilibrium.

\subsubsection{Hydrostatic equilibrium}
\label{sec:numerical_hydrostatic}

We consider the isentropic hydrostatic atmosphere \eqref{eq:isentropic_hydrostatic_equilibrium_sol} with $\hat{H}_0 = 1$ and $s_0 = 1$, thus it is given by
\begin{equation}
    \label{eq:init_hydro}
    \rho_\text{hyd}(x) = \left( 1 - \frac{\gamma - 1}{\gamma}\varphi(x)\right)^{1/(\gamma-1)},
    \quad
    p_\text{hyd}(x) = \rho(x)^\gamma,
    \quad
    u = 0.
\end{equation}
We use the two gravitational potentials $\varphi_1$ and $\varphi_2$.
The computational domain $\Omega = [0,1]$ is discretized in $50$ cells, thus constitutes a very coarse mesh.
Exact boundary conditions are prescribed for $\varphi_1$,
and periodic boundary conditions are prescribed for $\varphi_2$.
In \cref{tab:WB_hyro_linear,tab:WB_hydro_sine}, the~$L^2$ errors at final time~$t_f=1$
are compared to a non-well-balanced HLL solver for $\varphi_1$ and $\varphi_2$ respectively.
All well-balanced schemes are able to preserve the hydrostatic equilibrium at machine precision,
whereas the HLL scheme, as expected, has a quite large error.
Moreover, due to the non-trivial contribution of $\varphi_2$
the error using the HLL scheme is much higher.
This underlines the necessity of well-balanced schemes.

\begin{table}[tb]
    \renewcommand{\arraystretch}{1.25}
    \centering
    \caption{Well-balanced test case hydrostatic atmosphere with $\varphi(x) = (x-0.5)^2 / 2$.}
    \label{tab:WB_hyro_linear}
    \begin{tabular}{ccccc}
        \toprule
               & HLL                  & $\mathbb{P}_0$        & $\mathbb{P}_1^{\WB}$  & $\mathbb{P}_2^{\WB}$  \\
        \cmidrule(lr){1-5}
        $\rho$ & $6.390\cdot 10^{-3}$ & $2.658\cdot 10^{-15}$ & $2.727\cdot 10^{-15}$ & $2.436\cdot 10^{-15}$ \\
        $q$    & $5.002\cdot10^{-4}$  & $1.154\cdot 10^{-15}$ & $1.498\cdot 10^{-15}$ & $1.691\cdot 10^{-15}$ \\
        $E$    & $2.289\cdot10^{-3}$  & $7.640\cdot 10^{-15}$ & $7.640\cdot 10^{-15}$ & $1.020\cdot 10^{-14}$ \\
        \bottomrule
    \end{tabular}
\end{table}

\begin{table}[tb]
    \renewcommand{\arraystretch}{1.25}
    \centering
    \caption{Well-balanced test case hydrostatic atmosphere with $\varphi(x) = \sin(x)$.}
    \label{tab:WB_hydro_sine}
    \begin{tabular}{ccccc}
        \toprule
               & HLL                  & $\mathbb{P}_0$        & $\mathbb{P}_1^{\WB}$  & $\mathbb{P}_2^{\WB}$  \\
        \cmidrule(lr){1-5}
        $\rho$ & $1.109\cdot 10^{-1}$ & $3.819\cdot 10^{-15}$ & $2.620\cdot 10^{-15}$ & $4.255\cdot 10^{-15}$ \\
        $q$    & $1.898\cdot10^{-2}$  & $1.439\cdot 10^{-15}$ & $1.477\cdot 10^{-15}$ & $1.522\cdot 10^{-15}$ \\
        $E$    & $7.170\cdot10^{-2}$  & $7.536\cdot 10^{-15}$ & $6.822\cdot 10^{-15}$ & $7.270\cdot 10^{-15}$ \\
        \bottomrule
    \end{tabular}
\end{table}

\subsubsection{Moving equilibrium}
\label{sec:numerical_moving}

Next, a moving equilibrium is computed,
characterized by the triplet $q_0 = 1$, $s_0 = 0$ and $H_0 = 5$.
The gravitational potentials $\varphi_1$ and $\varphi_2$ are considered.
For each $\varphi$, the initial conditions for $\rho$ and $E$
are obtained by solving the nonlinear system \eqref{eq:equilibrium} employing Newton's method.
The computational domain is given by $\Omega = [0,1]$,
and is partitioned in $50$ cells.
Homogeneous exact and periodic boundary conditions
are prescribed for $\varphi_1$ and $\varphi_2$, respectively.
The $L^2$ errors at time $t=1$ are given in
\cref{tab:WB_moving_linear,tab:WB_moving_sine} for $\varphi_1$ and $\varphi_2$ respectively.
Similarly to the simulation of the hydrostatic atmosphere from \cref{sec:numerical_hydrostatic},
all well-balanced schemes are able to preserve the moving equilibrium up to machine precision,
whereas the HLL scheme yields quite large errors.
This verifies and illustrates the ability of the new well-balanced solvers
to capture moving equilibria up to machine precision.

\begin{table}[tb]
    \renewcommand{\arraystretch}{1.25}
    \centering
    \caption{Well-balanced test case moving equilibrium with $\varphi(x) = (x-0.5)^2 / 2$.}
    \label{tab:WB_moving_linear}
    \begin{tabular}{ccccc}
        \toprule
               & HLL                  & $\mathbb{P}_0$        & $\mathbb{P}_1^{\WB}$   & $\mathbb{P}_2^{\WB}$   \\
        \cmidrule(lr){1-5}
        $\rho$ & $5.569\cdot 10^{-3}$ & $3.204\cdot 10^{-15}$ & $ 2.897\cdot 10^{-15}$ & $ 2.907\cdot 10^{-15}$ \\
        $q$    & $1.129\cdot 10^{-3}$ & $1.916\cdot 10^{-15}$ & $ 1.727\cdot 10^{-15}$ & $ 1.803\cdot 10^{-15}$ \\
        $E$    & $4.188\cdot 10^{-3}$ & $7.742\cdot 10^{-15}$ & $ 8.970\cdot 10^{-15}$ & $ 6.764\cdot 10^{-15}$ \\
        \bottomrule
    \end{tabular}
\end{table}

\begin{table}[tb]
    \renewcommand{\arraystretch}{1.25}
    \centering
    \caption{Well-balanced test case moving equilibrium with $\varphi(x) = \sin(x)$.}
    \label{tab:WB_moving_sine}
    \begin{tabular}{ccccc}
        \toprule
               & HLL                    & $\mathbb{P}_0$         & $\mathbb{P}_1^{\WB}$   & $\mathbb{P}_2^{\WB}$   \\
        \cmidrule(lr){1-5}
        $\rho$ & $ 8.9233\cdot 10^{-2}$ & $ 2.627\cdot 10^{-15}$ & $ 3.607\cdot 10^{-15}$ & $ 4.246\cdot 10^{-15}$ \\
        $q  $  & $ 2.0153\cdot 10^{-2}$ & $ 1.877\cdot 10^{-15}$ & $ 1.727\cdot 10^{-15}$ & $ 1.922\cdot 10^{-15}$ \\
        $E  $  & $ 6.7186\cdot 10^{-2}$ & $ 9.930\cdot 10^{-15}$ & $ 9.357\cdot 10^{-15}$ & $ 9.607\cdot 10^{-15}$ \\
        \bottomrule
    \end{tabular}
\end{table}

\subsection{Perturbation of equilibrium states}
\label{sec:numerical_WB_pert}

As the main motivation behind the construction of well-balanced schemes lies
in the resolution of small perturbations around equilibria,
the next two sets of test cases concern the performance of the novel
fully well-balanced solver on such flows.
In particular, this allows the use of coarse meshes,
as the solution is not polluted by background errors stemming
from the resolution of the equilibrium state,
which would then require a substantial grid refinement to
reduce the truncation error and make the perturbation visible.

In this section, we consider the quadratic potential $\varphi_1$,
and the computational domain is given by $\Omega = [0,1]$,
discretized using $50$ cells, with inhomogeneous Dirichlet boundary conditions
corresponding to the exact, unperturbed steady solution.
For all cases, we set
$C_\theta = 1$ for the $\mathbb{P}_1^\WB$ scheme
and $C_\theta = 0.15$ for the $\mathbb{P}_2^\WB$ scheme.

\subsubsection{Perturbations of a hydrostatic equilibrium}
\label{sec:numerical_hydrostatic_perturbed}

\begin{figure}[htbp]
    \centering
    \includegraphics[width=0.95\textwidth]{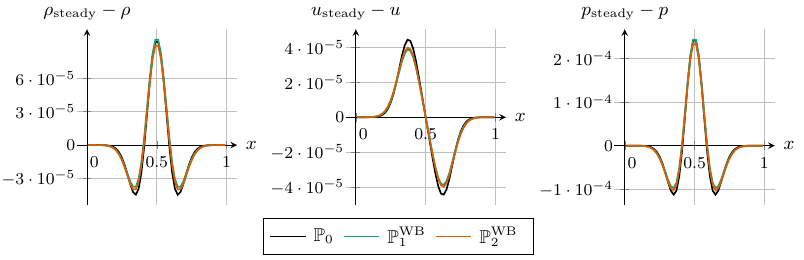}
    \includegraphics[width=0.95\textwidth]{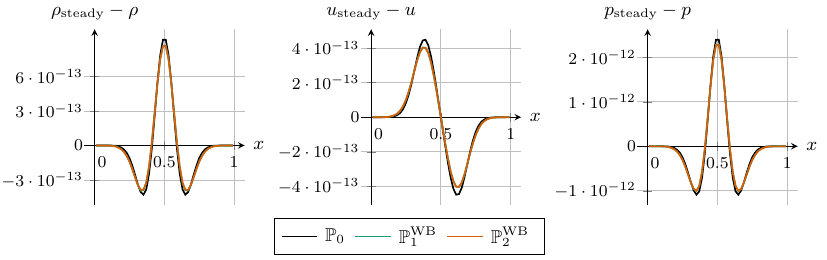}
    \caption{%
    Perturbation of a hydrostatic equilibrium.
    Top panels: perturbation with initial amplitude $A=10^{-4}$.
    Bottom: perturbation with initial amplitude $A=10^{-12}$.%
    }
    \label{fig:hydro_perturbed}
\end{figure}

As a first example, we consider a perturbation around an isentropic
atmosphere~\eqref{eq:init_hydro} taken e.g. from \cite{ChaKli2015}.
The initial condition is given by
\begin{equation*}
    \rho(x,0) = \rho_\text{hyd}(x),
    \
    u(x,0) = 0,
    \
    p(x,0) = p_\text{hyd}(x) + A \exp\left(-100 \left( x - \frac{1}{2}\right)^2\right),
\end{equation*}
where we set the amplitudes of the perturbations to be $A = 10^{-4}$ and $10^{-12}$,
where the latter is of the order of the truncation error of the well-balanced solver.
The difference of the numerical solution against the equilibrium
$W_\text{steady} - W_{\text{num}}$ obtained with the novel well-balanced
solvers at final time $t = 0.075$ are depicted in \cref{fig:hydro_perturbed}.
The top panels contain the results with an initial perturbation of amplitude $10^{-4}$,
while the bottom ones display the results with the amplitude reduced to $10^{-12}$.
In both cases, we observe that the perturbations are clearly visible and
there are no spurious artifacts from the background atmosphere.

\subsubsection{Perturbations of a moving equilibrium}
\label{sec:numerical_moving_perturbed}

\begin{figure}[htbp]
    \centering
    \includegraphics[width=0.95\textwidth]{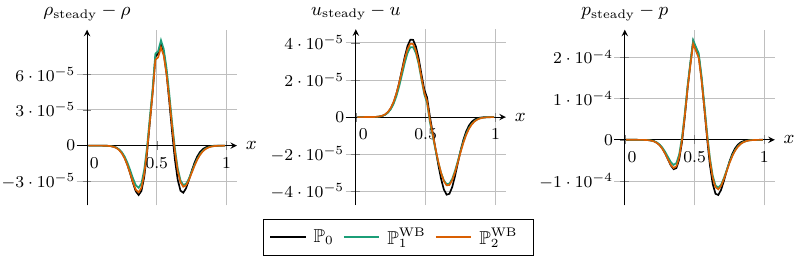}
    \includegraphics[width=0.95\textwidth]{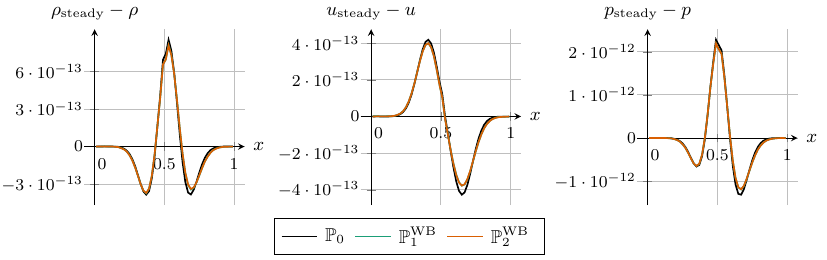}
    \caption{%
    Perturbation of a moving equilibrium.
    Top panels: perturbation with initial amplitude $A=10^{-4}$.
    Bottom: perturbation with initial amplitude $A=10^{-12}$.%
    }
    \label{fig:moving_perturbed}
\end{figure}

Next, a perturbation of the moving equilibrium is considered.
The initial condition is given by
\begin{equation*}
    \rho(x,0) = \rho_{0}(x),
    \
    q(x,0) = q_0,
    \
    p(x,0) = p_0(x) + A \exp\left(-100 \left( x - \frac{1}{2}\right)^2\right),
\end{equation*}
where the equilibrium states are obtained from the triplet $(q_0,s_0,H_0)$
already introduced to build the steady solution in \cref{sec:numerical_moving}.
The numerical results are given in \cref{fig:moving_perturbed} at the final time $t = 0.075$,
arranged similarly to \cref{fig:hydro_perturbed}.
The perturbation around the moving equilibrium is once again well-captured
and no spurious errors are introduced from the background equilibrium.
In comparison to the hydrostatic atmosphere case from \cref{fig:hydro_perturbed},
the perturbations are skewed towards the right and travel faster.
This is due to a non-zero positive background velocity.

\begin{figure}[htbp]
    \centering
    \includegraphics[width=0.95\textwidth]{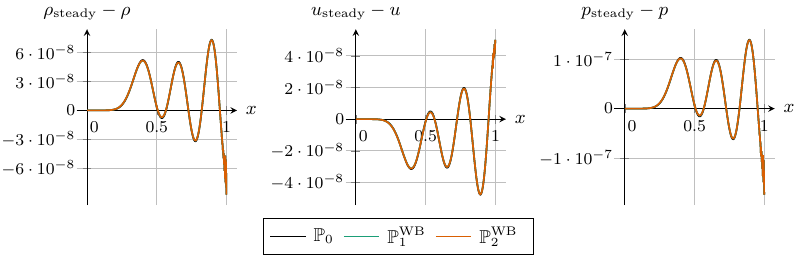}
    \caption{Perturbation of moving equilibrium from the boundary.}
    \label{fig:moving_perturbed_boundary}
\end{figure}

Another interesting test case concerns a moving steady solution perturbed by
a wave stemming from a time-dependent boundary condition.
We consider again the moving equilibrium given by $(q_0, s_0, H_0)$,
with the perturbation applied as a right boundary condition acting only on the velocity,
which is given at $x=1$ and $t>0$ by
\begin{equation*}
    u(1,t) = 10^{-8} \sin(8 \pi t)
\end{equation*}
instead of by the unperturbed equilibrium.
The computational domain is given by $[0,1]$ using $512$ cells
and the simulation is stopped shortly before the perturbation reaches the opposite boundary,
at $t = 0.72$.
Similar set-ups can be found in \cite{FucMcmMisRisWaa2010},
motivated by the study of wave propagation in stellar atmospheres.
The difference between numerical results and the moving equilibrium
are given in \cref{fig:moving_perturbed_boundary},
where we have set $\Lambda = 5$ to add some diffusion to the scheme.
Indeed, in this case, the solution is rapidly changing,
due to the high-frequency oscillation introduced from the boundary.
This can be compared to a solution with strong gradients, for which it is well-known that additional diffusion is required to avoid spurious oscillations.
Therefore, in this test case,
we have chosen a larger value of $\Lambda$ to stabilize the scheme.
Despite this extra diffusion,
we observe that the perturbation remains well-resolved.
This further verifies the ability of the numerical scheme to resolve
small perturbations around hydrostatic and moving equilibria on coarse grids.

\subsection{Riemann problems}
\label{sec:numerical_RP}

As a final series of test cases, we consider four Riemann problems, three classical ones for the Euler equations without a source term,
and one in presence of a gravitational field, far from an equilibrium.
Although the first three Riemann problems do not involve a gravitational field,
they serve as validation of the consistency properties of our new well-balanced scheme towards the homogeneous Euler equations in the limit of vanishing gravitational influence.
Indeed, if the source term vanishes,
our scheme should reduce to a classical HLL scheme,
and be able to yield good approximations of the exact solution.

In all these cases, the space domain is $\Omega = [0,1]$,
and we prescribe homogeneous Neumann boundary conditions.
Moreover, the initial condition takes the form of a Riemann problem, i.e.,
\begin{equation*}
    W(x,0) = \begin{cases}
        W_L & \text{if } x < x_0,    \\
        W_R & \text{if } x \geq x_0,
    \end{cases}
\end{equation*}
with the jump position $x_0 = 0.5$,
and where the left and right states are given differently for each problem.
In each case, we take a grid consisting of $75$ cells.

\subsubsection{Sod problem}

\begin{figure}[tb]
    \centering
    \includegraphics[width=0.8\textwidth]{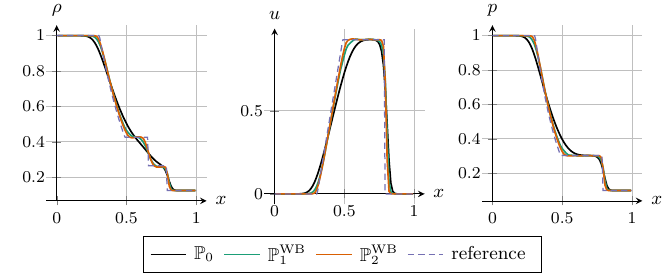}
    \caption{%
        Riemann problem without a gravitational field influence:
        Sod problem displayed in primitive variables $\rho$, $u$ and $p$,
        at time $t = 0.1644$ using $75$ cells.%
    }
    \label{fig:RP_sod}
\end{figure}

We first consider the standard Sod problem, which is a well-known benchmark for Riemann solvers.
The initial condition is given by
\begin{equation*}
    \rho_L = 1, \quad \rho_R = 0.125, \quad u_L = u_R = 0, \quad p_L = 1, \quad p_R = 0.1.
\end{equation*}
The gravitational potential is given by $\varphi=0$.
The results are given in \cref{fig:RP_sod},
in primitive variables $\rho$, $u$ and $p$ at time $t = 0.1644$.
The first order scheme is diffusive, in particular over the contact and shock waves,
whereas increasing the order of the new well-balanced scheme increases
the resolution of all waves.
As is typical for high-order schemes,
small oscillations can be observed near the shock wave.
However, all schemes are able to correctly determine the shock position and amplitude.
The numerical results are compared against an exact reference solution obtained with an exact Riemann solver \cite{Toro2009}.

\subsubsection{Double rarefaction}

\begin{figure}[tb]
    \centering
    \includegraphics[width=0.8\textwidth]{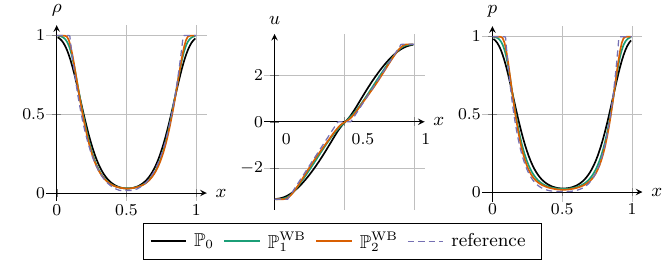}
    \caption{%
        Riemann problem without a gravitational field influence:
        double rarefaction problem displayed in primitive variables $\rho$, $u$ and $p$,
        at time $t = 0.09$ using $75$ cells.%
    }
    \label{fig:RP_vacuum}
\end{figure}

We then turn to a double rarefaction, whose initial data is
\begin{equation*}
    \rho_L = \rho_R = 1, \quad u_L = -\frac{10}3, u_R = \frac{10}3, \quad p_L = p_R = 1.
\end{equation*}
This is a challenging problem, as it will lead to a near-vacuum state in the middle of the domain.
We still take the given gravitational potential $\varphi=0$.
The results are displayed in \cref{fig:RP_vacuum} at time $t = 0.09$.
Despite the fact that $\rho$ and $p$ are very close to zero in the middle of the domain,
no negative values of $\rho$ and $p$ are observed.
The numerical results are compared against an exact reference solution obtained with an exact Riemann solver \cite{Toro2009}.

\subsubsection{Stationary shock wave}
\begin{figure}[tb]
    \centering
    \includegraphics[width=0.8\textwidth]{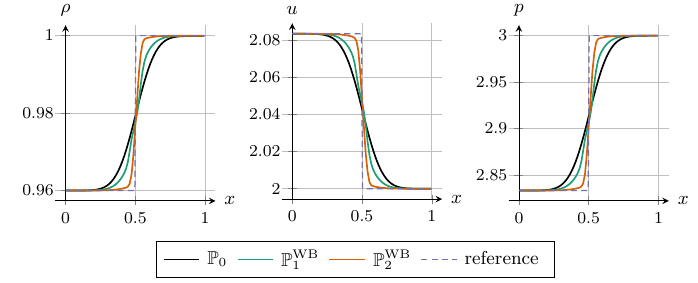}
    \caption{%
        Riemann problem without a gravitational field influence:
        stationary shock wave displayed in primitive variables $\rho$, $u$ and $p$,
        at time $t = 0.25$ using $75$ cells.%
    }
    \label{fig:RP_shock}
\end{figure}

The initial condition of the third Riemann problem is designed to have a single stationary shock wave which is given by
\begin{equation*}
    \begin{matrix}
        \rho_L = 24/25, \quad & u_L = 25/12, \quad & p_L = 17/6, \\
        \rho_R = 1            & u_R = 2            & p_R = 3.
    \end{matrix}
\end{equation*}
The purpose of this test is to assess the behavior of the numerical scheme
in the presence of such a stationary discontinuity in the absence of gravitational influence, i.e. $\varphi=0$.
The results are plotted in
\cref{fig:RP_shock}
at time $t = 0.25$.
This time, we take $C_\theta = 3$ for both the $\mathbb{P}_1^\WB$ and $\mathbb{P}_2^\WB$ schemes.
For all schemes, we observe good agreement with the exact solution.
However, as expected, this stationary solution is not exactly captured by our well-balanced scheme.
This is not in contradiction with the above derived method, since it is designed to be exact on smooth moving equilibria given by \eqref{eq:equilibrium}.
While the constraint on the momentum $q_0$ is satisfied, the enthalpy constraint $H_0$ is not.

\subsubsection{Riemann problem in a gravitational field}

\begin{figure}[htbp]
    \centering
    \includegraphics[width=0.8\textwidth]{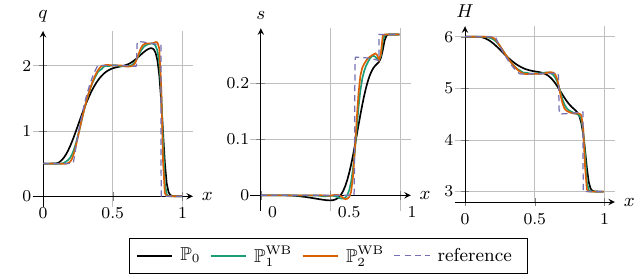}
    \caption{%
        Riemann problem under a gravitational field influence:
        modified Sod problem displayed in steady variables $q$, $s$ and $H$,
        at time $t = 0.2$ using $75$ cells.%
    }
    \label{fig:RP_sod_source}
\end{figure}

The last Riemann problem is a Sod-like problem in the presence of a gravitational source term,
given by $\varphi_1$.
The initial data is given in steady variables $q$, $s$ and $H$, as follows:
\begin{equation*}
    q_L = 0.5, \; q_R = 0,
    \quad
    s_L = 0, \; s_R = -\log(0.75),
    \quad
    H_L = 6, \; H_R = 3.
\end{equation*}
The results are displayed in \cref{fig:RP_sod_source} at time $t = 0.2$.
This time, the equilibrium variables $q$, $s$ and $H$ are displayed
corresponding to the variables, in which the initial data is given.
We observe that the solution generated by the new well-balanced scheme accurately captures the arising shock and rarefaction waves.
The numerical results are compared against a reference solution
obtained with the classical HLL scheme with a centered discretization of the source terms on $2000$ cells.

\subsection{Two-dimensional extension}
\label{sec:2D}

We close this series of test cases by providing
a two-dimensional (2D) extension of the well-balanced scheme.
Since it is designed to exactly capture
isentropic moving equilibria in one dimension,
the resulting scheme will not be well-balanced
for generic 2D equilibria.
Rather, it will exactly preserve one-dimensional (1D)
moving equilibria in each direction;
it will also exactly capture all isentropic hydrostatic equilibria.
Moreover, according to e.g. \cite{PerShu1996},
2D schemes can be recast as convex combinations of 1D schemes.
Therefore, all robustness properties of the 1D scheme
will be inherited by the 2D scheme.
Namely, the 2D scheme will also be entropy-stable
and positivity-preserving.

Based on this remark,
we derive a straightforward 2D extension of
the scheme given by \eqref{eq:scheme_from_ARS}.
The 2D version of system \eqref{eq:EulerG} reads
\begin{equation}
    \label{eq:EulerG_fully_2D}
    \begin{dcases}
        \pd_t \rho
        + \nabla \cdot (\rho \boldsymbol u) =
        0 ,                    \\
        \pd_t (\rho \boldsymbol u)
        + \nabla \cdot (\rho \boldsymbol u \otimes \boldsymbol u)
        + \nabla p =
        - \rho \nabla \varphi, \\
        \pd_t E
        + \nabla \cdot ((E + p) \boldsymbol u) =
        - \rho \boldsymbol u \cdot \nabla \varphi,
    \end{dcases}
\end{equation}
where
\begin{equation}
    \label{eq:2D_pressure_law}
    p = (\gamma - 1) \left( E - \frac 1 2 \rho \|\boldsymbol u\|^2 \right).
\end{equation}
To approximate solutions to this system,
we consider a Cartesian mesh with space steps $\dx$ and $\dy$,
and with cells $(x_{\imh,j}, x_{\iph,j}) \times (y_{i,\jmh}, y_{i,\jph})$.
The approximate solution $W_{i,j}^n$ of
\eqref{eq:EulerG_fully_2D} is defined by
the cell average
\begin{equation*}
    \label{eq:cell_average_W_2D}
    W_{i,j}^n \simeq
    \frac 1 \dx
    \frac 1 \dy
    \int_{x_{\imh,j}}^{x_{\iph,j}}
    \int_{y_{i,\jmh}}^{y_{i,\jph}}
    W(x,y,t^n) \, dy \, dx.
\end{equation*}
We denote the unit normal vectors in the $x$ and $y$ directions by
\smash{$\normalx = \begin{pmatrix}1 & 0 \end{pmatrix}^T$}
and \smash{$\normaly = \begin{pmatrix}0 & 1 \end{pmatrix}^T$}.
Then, this cell average is evolved in time according to
\begin{equation*}
    \label{eq:2D_scheme_from_ARS}
    \begin{aligned}
        W_{i,j}^{n+1} = W_{i,j}^n + \frac \dt \dx \Big[
         & \lambda_{\iph,j}
        \Big(W_L^* \left( W_{i,j}^{n}, W_{i+1,j}^n, \varphi_{i,j}, \varphi_{i+1,j}, \normalx \right) - W_{i,j}^n\Big) \\
            +
         & \lambda_{\imh,j}
            \Big(W_R^* \left( W_{i-1,j}^{n}, W_{i,j}^n, \varphi_{i-1,j}, \varphi_{i,j}, \normalx \right) - W_{i,j}^n\Big)
        \Big]                                                                                                         \\
        + \frac \dt \dy \Big[
         & \lambda_{i,\jph}
        \Big(W_L^* \left( W_{i,j}^{n}, W_{i,j+1}^n, \varphi_{i,j}, \varphi_{i,j+1}, \normaly \right) - W_{i,j}^n\Big) \\
            +
         & \lambda_{i,\jmh}
            \Big(W_R^* \left( W_{i,j-1}^{n}, W_{i,j}^n, \varphi_{i,j-1}, \varphi_{i,j}, \normaly \right) - W_{i,j}^n\Big)
            \Big],
    \end{aligned}
\end{equation*}
where the intermediate states
\begin{equation}
    \label{eq:intermediate_states_2D}
    W_L^*\big(W_L, W_R, \varphi_L, \varphi_R, \boldsymbol n\big)
    \text{\quad and \quad}
    W_R^*\big(W_L, W_R, \varphi_L, \varphi_R, \boldsymbol n\big)
\end{equation}
are defined by a forthcoming 2D extension of the 1D
\ARS \eqref{eq:ARS},
according to some given normal
unit vector $\boldsymbol n$.
The approximate wave speeds and time step
are defined by straightforward extensions of the 1D case.

To derive the intermediate states \eqref{eq:intermediate_states_2D},
we consider the 2D system
written in an arbitrary direction $\boldsymbol n$, see e.g. \cite{Toro2009}.
This leads to the augmented one-dimensional system
\begin{equation}
    \label{eq:EulerG_2D}
    \begin{dcases}
        \pd_t \rho  + \pd_x (\rho u) = 0 ,                             \\
        \pd_t (\rho u)  + \pd_x (\rho u^2 + p) = - \rho \pd_x \varphi, \\
        \pd_t (\rho v)  + \pd_x (\rho u v) = 0,                        \\
        \pd_t E  + \pd_x ((E + p) u) = - \rho u \pd_x \varphi,
    \end{dcases}
\end{equation}
where $u = \boldsymbol u \cdot \boldsymbol n$
is the normal velocity and
$v = \boldsymbol u \cdot \boldsymbol n^\perp$
is the tangential velocity.
For simplicity, the variable $x$ denotes the normal direction.
The steady solutions of \eqref{eq:EulerG_2D}
are given by \eqref{eq:equilibrium},
with the additional relation $v \eqqcolon v_0$
describing a constant tangential velocity.
That is to say, the steady solutions are characterized by
\begin{equation*}
    \rho u = q_0,
    \quad
    -\log \left( \frac p {\rho^\gamma} \right) = s_0,
    \quad
    \frac{E + p}{\rho} + \varphi = H_0,
    \quad
    v = v_0,
\end{equation*}
with $q_0$, $s_0$, $H_0$ and $v_0$ constant values.
To derive a fully well-balanced \ARS for~\eqref{eq:EulerG_2D},
we first directly reuse the 1D intermediate states
for $\rho$, $q = \rho u$ and $E$,
taking care to use the 2D pressure law \eqref{eq:2D_pressure_law}.
Second, we derive a new intermediate state for $\rho v$.
These intermediate states should be consistent and well-balanced,
i.e., they should exactly preserve constant tangential velocities.
Straightforward computations lead to
\begin{equation*}
    (\rho v)_L^* = \widehat{(\rho v)} - \delta(\rho v),
    \qquad
    (\rho v)_R^* = \widehat{(\rho v)} + \delta(\rho v),
\end{equation*}
with
\begin{equation*}
    \widehat{(\rho v)} = (\rho v)_\HLL
    \text{\quad and \quad}
    \delta(\rho v) = \delta \rho \, \frac{(\rho v)_\HLL}{\rho_\HLL}.
\end{equation*}

Thus, the scheme \eqref{eq:2D_scheme_from_ARS}
has been fully defined.
To properly describe its properties, we first state the following definition
of grid-aligned equilibria.

\begin{definition}
    \label{def:2D_grid_aligned_equilibria}
    Consider the 2D Euler system \eqref{eq:EulerG_fully_2D}
    with $\boldsymbol{u} = (u_1, u_2)^T$.
    We call grid-aligned equilibrium any solution such that either
    \begin{equation*}
        \partial_y \varphi = 0,
        \quad
        \rho u_1 = q_0,
        \quad
        u_2 = v_0,
        \quad
        -\log \left( \frac p {\rho^\gamma} \right) = s_0,
        \quad
        \frac{E + p}{\rho} + \varphi = H_0,
    \end{equation*}
    or
    \begin{equation*}
        \partial_x \varphi = 0,
        \quad
        \rho u_2 = q_0,
        \quad
        u_1 = v_0,
        \quad
        -\log \left( \frac p {\rho^\gamma} \right) = s_0,
        \quad
        \frac{E + p}{\rho} + \varphi = H_0.
    \end{equation*}
\end{definition}

Equipped with this definition,
the properties of the 2D scheme are summarized in the following result.
Its proof is a straightforward extension of the 1D case,
and of the recasting of 2D schemes as convex combinations of 1D schemes,
see \cite{PerShu1996}.

\begin{theorem}
    \label{thm:2D_scheme}
    The scheme \eqref{eq:2D_scheme_from_ARS} is entropy-stable,
    positivity-preserving and
    preserves all grid-aligned moving equilibria.
\end{theorem}

To test the properties outlined in \cref{thm:2D_scheme},
we consider an exact solution in \cref{sec:GravGresho},
two equilibrium states in \cref{sec:2D_WB_stationary,sec:2D_steady_perturbed},
respectively unperturbed and perturbed,
and lastly an implosion in \cref{sec:implosion}.

\subsubsection{Stationary vortex in a gravitational potential}
\label{sec:GravGresho}
In this test case, we consider a vortex in a gravitational field from \cite{Thomann2020}.
It is a continuous stationary solution of the Euler equations with gravity, based on an isothermal hydrostatic equilibrium, with a closed-form expression.
The purpose of the test is to show the convergence of the first order scheme in a two-dimensional set-up.
Note that the initial condition is only continuous and not continuously differentiable.
Thus, even for higher order schemes, a full order of convergence cannot be expected and we therefore only show the rates of the first order scheme.
Moreover, since the underlying equilibrium is not isentropic, this test checks whether the scheme has good results on non-isentropic equilibria.

The initial condition is given in polar coordinates, with respect to the radius $r = \sqrt{x-x_0)^2 + (y-y_0)^2}$ and the angle $\theta = \arctan((y-y_0)/(x-x_0))$.
The initial condition satisfies
\begin{equation*}
    \partial_r p = \frac{\rho u_\theta^2}{r} - \rho \partial_r \phi.
\end{equation*}
Therein, $u_\theta$ denotes the angular velocity of the vortex, is given by
\begin{equation*}
    u_\theta = \begin{cases}
        5 r     & \text{ if } r \leq \frac{1}{5},               \\
        2 - 5 r & \text{ if } \frac{1}{5} < r \leq \frac{2}{5}, \\
        0       & \text{ if } r > \frac{2}{5}.
    \end{cases}
\end{equation*}
The density is given by a isothermal hydrostatic atmosphere $\rho_\text{hyd}$, and the pressure $p = p_\text{hyd} + p_\text{vort}$ is a perturbation of the hydrostatic one.
The hydrostatic density and pressure are characterized by
\begin{equation*}
    \rho_\text{hyd}(r) = \exp\left(-\frac{\varphi(r)}{RT}\right),
    \text{\qquad and \qquad}
    p_\text{hyd}(r) = RT \rho_\text{hyd}(r),
\end{equation*}
with $R$ and $T$ are constant values, defined such that $RT = 4$.
Then, $p_\text{vort}$ has to satisfy
\begin{equation*}
    p_\text{vort} = \int_0^r \frac{1}{s} \rho_\text{hyd}(s) u_\theta^2(s) ds.
\end{equation*}
To compute density and pressure, we define a piecewise continuous gravitational potential $\varphi$ as
\begin{equation*}
    \varphi(r) = \begin{cases}
        \frac{25}{2} r^2                                                                                                      & \text{ if } r \leq \frac{1}{5}               \\
        \frac{1}{2} - \ln(\frac{1}{5}) + \ln(r)                                                                               & \text{ if } \frac{1}{5} < r \leq \frac{2}{5} \\
        \ln(2) - \frac{1}{2} \frac{r_c}{r_c - 2/5}  + \frac{5}{2} \frac{r_c}{r_c - 2/5}r - \frac{5}{4}\frac{1}{r_c - 2/5} r^2 & \text{ if } \frac{2}{5} < r \leq r_c         \\
        \ln(2) - \frac{1}{2} \frac{r_c}{r_c - 2/5} + \frac{5}{4}\frac{r_c^2}{r_c - 2/5}                                       & \text{ if } r > r_c.
    \end{cases}
\end{equation*}
This yields a pressure $p_\text{vort}$ defined in a piecewise fashion by
\begin{equation*}
    p_\text{vort} = RT\begin{cases}
        p_1(r)                              & \text{ if } r \leq \frac{1}{5}               \\
        p_1(\frac{1}{5}) + p_2(r)           & \text{ if } \frac{1}{5} < r \leq \frac{2}{5} \\
        p_1(\frac{1}{5}) + p_2(\frac{2}{5}) & \text{ if } r > \frac{2}{5},
    \end{cases}
\end{equation*}
with
\begin{equation*}
    \begin{aligned}
        p_1(r) & = 1 - \exp\left(-\frac{25r^2}{2RT}\right), \\
        p_2(r) & =
        p_{2, 1}
        \biggl(
        p_{2,2}+
        r ^{\frac{-1}{RT}}
        \Bigl( - 2 +
        10 r (1 - 2 RT)
        + RT (6 - 4 RT)
        + \frac{25r^2}{2} (RT - 1)\Bigr)
        \biggr),
    \end{aligned}
\end{equation*}
where we have set
\begin{equation*}
    p_{2, 1} = \frac{\exp\left(-\frac{1/2 + \ln(5)}{RT}\right)}{(RT - 1) (RT - 1/2)}
    \text{\quad and \quad}
    p_{2, 2} = \exp\left(\frac{\ln(5)}{RT}\right) \left(RT \left(4 RT - \frac{5}{2}\right) + \frac{1}{2}\right).
\end{equation*}

In the simulation, we set $r_c = \frac{3}{5}$ and the computational domain is $[0,1]^2$, with periodic boundary conditions.
Note that the vortex is constructed such that it is constant near the boundary.
In \cref{tab:errors_p} the errors in the $L^2$ and the $L^\infty$ norm are given.
As expected, under mesh refinement, the experimental order of convergence approaches the theoretical first order of accuracy.
Similar results are obtained for the density and velocity and are thus omitted.

\begin{table}[tb]
    \renewcommand{\arraystretch}{1.25}
    \centering
    \caption{Stationary vortex in a gravitational potential from \cref{sec:GravGresho}: errors for the pressure $p$ in $L^2$ and $L^\infty$ norms.}
    \label{tab:errors_p}
    \begin{tabular}{cccccc}
        \toprule
               & \multicolumn{2}{c}{$L^2$ error} & \multicolumn{2}{c}{$L^\infty$ error}                                \\
        \cmidrule(lr){2-3} \cmidrule(lr){4-5}
        $N$    & Error                           & Order                                & Error               & Order  \\
        \cmidrule(lr){1-5}
        $64$   & $6.23\cdot 10^{-2}$             & ---                                  & $5.01\cdot 10^{-1}$ & ---    \\
        $128$  & $4.49\cdot 10^{-2}$             & $0.47$                               & $3.33\cdot 10^{-1}$ & $0.59$ \\
        $256$  & $2.79\cdot 10^{-2}$             & $0.69$                               & $1.86\cdot 10^{-1}$ & $0.84$ \\
        $512$  & $1.59\cdot 10^{-2}$             & $0.81$                               & $9.83\cdot 10^{-2}$ & $0.92$ \\
        $1024$ & $8.63\cdot 10^{-3}$             & $0.88$                               & $5.11\cdot 10^{-2}$ & $0.94$ \\
        \bottomrule
    \end{tabular}
\end{table}

\subsubsection{Well-balanced assessment: preservation of a stationary solution}
\label{sec:2D_WB_stationary}
The next test case regards the exact capture of a stationary state.
Since the scheme is not designed to preserve truly two-dimensional moving steady states in a gravitational field $\varphi = \varphi(x,y)$ due to the $\nabla \cdot \rho \boldsymbol u = 0$ condition, we consider a moving stationary state in the gravitational field $\varphi(y) = \sin(2\pi y)$.
This means that, on the one hand, in the $y$-direction, we have a one-dimensional moving equilibrium with constant momentum.
On the other hand, in the $x$-direction, no gravitational influence is present, and the solution is merely transported at constant velocity along the $x$-axis.
The moving stationary state $W_\text{steady}$ for a velocity field $\boldsymbol u = (u_1, u_2)^T$ reads
\begin{equation*}
    u_1 = 1, \quad \rho u_2 = 1, \quad s_0 = 0, \quad H_0 = 10,
\end{equation*}
in the context of \cref{def:2D_grid_aligned_equilibria}.

The simulation is performed on a coarse mesh made of $32^2$ cells on the computational domain $[0,1]^2$, with periodic boundary conditions, and up to a final time $t_f =1$.
In \cref{tab:errors_rho_p}, we compare the results for the density and the pressure obtained with our new well-balanced scheme, against the non-well-balanced HLL scheme.
As expected, we observe that the HLL scheme generates quite large errors, whereas the well-balanced scheme is able to preserve the steady state up to machine precision.
Similar results are obtained for the velocity and momentum, and are thus omitted.

\begin{table}[tb]
    \renewcommand{\arraystretch}{1.25}
    \centering
    \caption{Preservation of the moving stationary solution from \cref{sec:2D_WB_stationary}: relative errors for the density $\rho$ and pressure $p$ in $L^2$ and $L^\infty$ norms.}
    \label{tab:errors_rho_p}
    \begin{tabular}{ccccc}
        \toprule
         & \multicolumn{2}{c}{$\rho$}
         & \multicolumn{2}{c}{$p$}                           \\
        \cmidrule(lr){2-3} \cmidrule(lr){4-5}
         & $L^2$ error                & $L^\infty$ error
         & $L^2$ error                & $L^\infty$ error     \\
        \cmidrule(lr){1-5}
        HLL
         & $4.90\cdot 10^{-3}$        & $9.30\cdot 10^{-3}$
         & $2.59\cdot 10^{-2}$        & $3.55\cdot 10^{-2}$  \\
        WB
         & $3.42\cdot 10^{-16}$       & $9.06\cdot 10^{-16}$
         & $4.13\cdot 10^{-16}$       & $1.10\cdot 10^{-15}$ \\
        \bottomrule
    \end{tabular}
\end{table}

\subsubsection{Well-balanced assessment: perturbation of a stationary solution}
\label{sec:2D_steady_perturbed}
Analogously to the one-dimensional test cases, we study the evolution of a perturbation of the two-dimensional moving steady state $W_\text{steady}$ from \cref{sec:2D_WB_stationary}.
The perturbation is given by a Gaussian density pulse
\begin{equation*}
    \rho = \rho_\text{steady} + 10^{-5} e^{-50 r^2},
\end{equation*}
with $r^2 = (x-0.5)^2 + (y-0.5)^2$, centered at $(0.5,0.5)$, in the middle of the computational domain $[0,1]^2$.

\begin{figure}[!ht]
    \centering
    \begin{subfigure}{0.45\textwidth}
        \includegraphics[width=\textwidth]{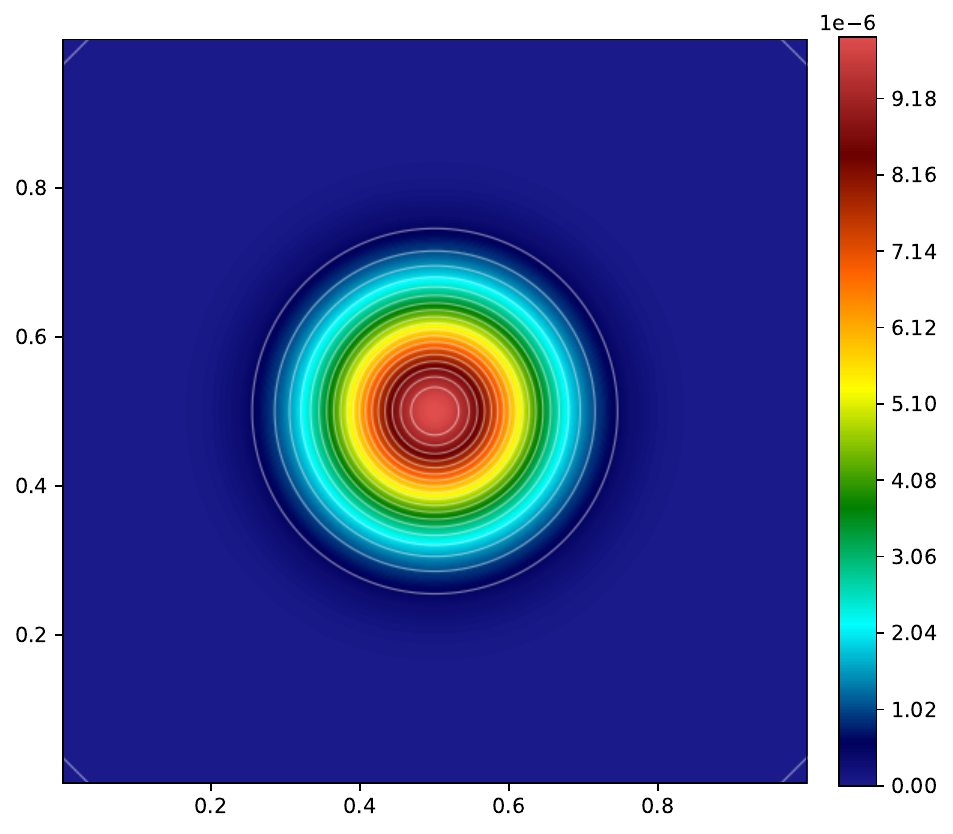}
        \caption{$t=0$}
        \label{fig:2D_steady_perturbed_1}
    \end{subfigure}
    \hfill
    \begin{subfigure}{0.45\textwidth}
        \includegraphics[width=\textwidth]{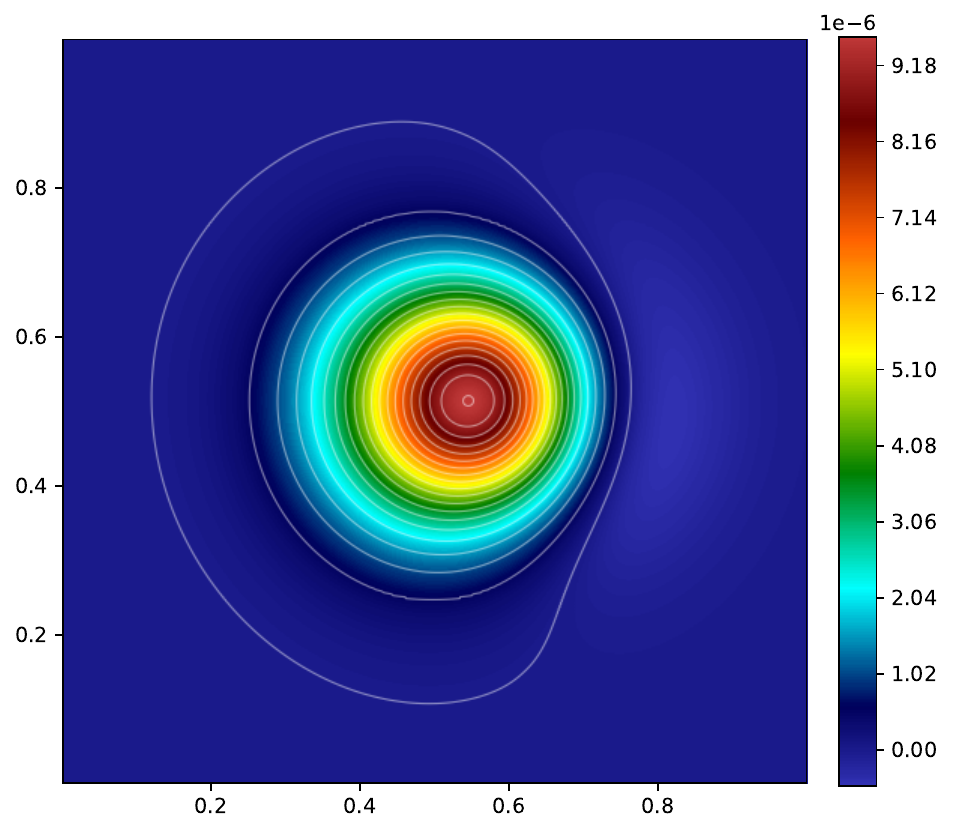}
        \caption{$t=t_f / 3$}
        \label{fig:2D_steady_perturbed_2}
    \end{subfigure}

    \bigskip

    \begin{subfigure}{0.45\textwidth}
        \includegraphics[width=\textwidth]{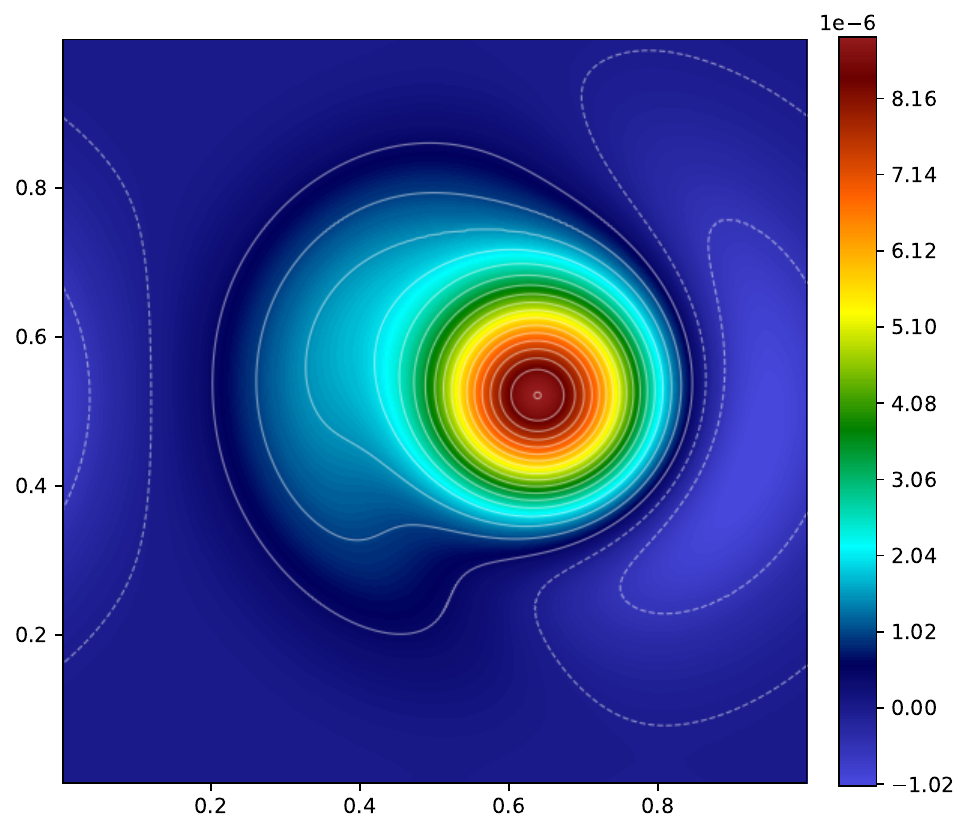}
        \caption{$t= 2 t_f / 3$}
        \label{fig:2D_steady_perturbed_3}
    \end{subfigure}
    \hfill
    \begin{subfigure}{0.45\textwidth}
        \includegraphics[width=\textwidth]{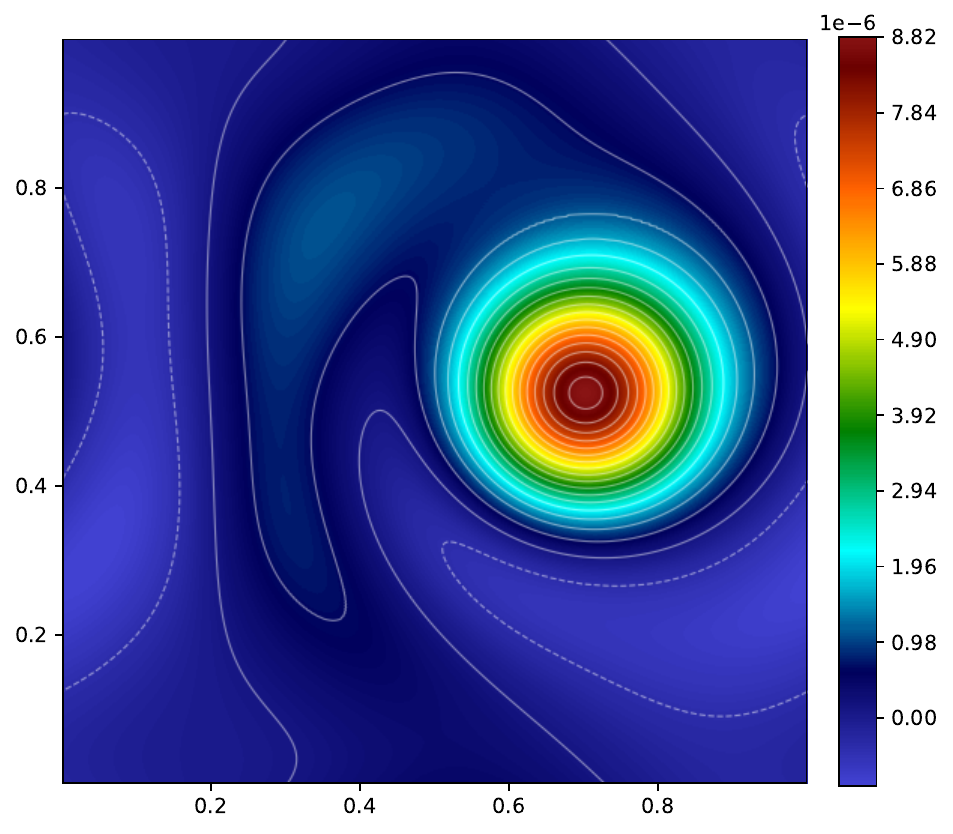}
        \caption{$t=t_f$}
        \label{fig:2D_steady_perturbed_4}
    \end{subfigure}

    \caption{Perturbed steady solution from \cref{sec:2D_steady_perturbed}: perturbation in density $ \rho - \rho_\text{steady}$ at different times, using a mesh made of $512^2$ cells. Note that the color scale is the same for all figures.}
    \label{fig:2D_steady_perturbed}
\end{figure}

In \cref{fig:2D_steady_perturbed}, the density perturbation $\rho - \rho_\text{steady}$ at different snapshots in time up to the final time $t_f = 0.2$ is displayed, for a mesh made of $512^2$ cells and with periodic boundary conditions.
Due to the underlying steady solution, the small density perturbation is transported along the $x$-axis without creating spurious artifacts from the underlying steady state.

\subsubsection{Implosion}
\label{sec:implosion}

The last two-dimensional test case under consideration is a Riemann Problem given by
a radial implosion under a Gaussian gravitational potential
\begin{equation*}
    \varphi(x,y) =
    - \exp \left( -50 \left( (x-0.5)^2 + (y-0.5)^2 \right) \right).
\end{equation*}
The initial conditions are given by
\begin{equation*}
    W(x,y,0) =
    \begin{dcases}
        \begin{pmatrix}
            1 & 0 & 0 & 1
        \end{pmatrix}^T
         & \text{ if } r < 0.25 \text{,} \\
        \begin{pmatrix}
            0.125 & 0 & 0 & 0.1
        \end{pmatrix}^T
         & \text{ otherwise,}
    \end{dcases}
\end{equation*}
depending on the radius $r = \sqrt{(x-0.5)^2 + (y-0.5)^2}$.

\begin{figure}[!ht]
    \centering
    \begin{subfigure}{0.45\textwidth}
        \includegraphics[width=\textwidth]{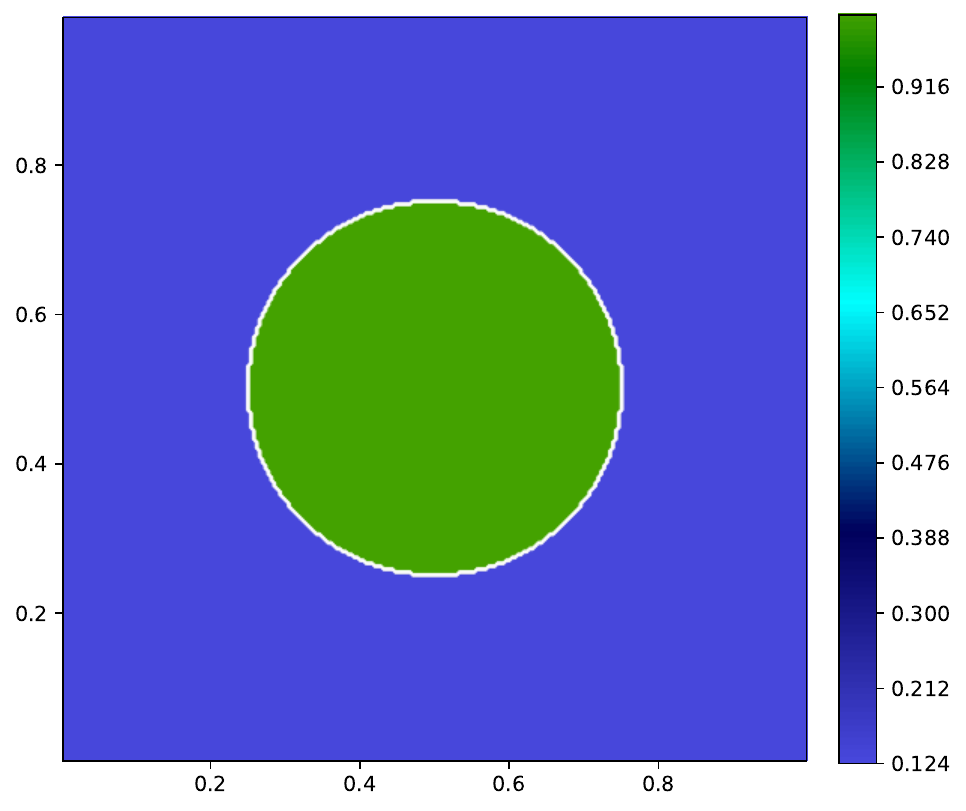}
        \caption{$t=0$}
        \label{fig:2D_implosion_1}
    \end{subfigure}
    \hfill
    \begin{subfigure}{0.45\textwidth}
        \includegraphics[width=\textwidth]{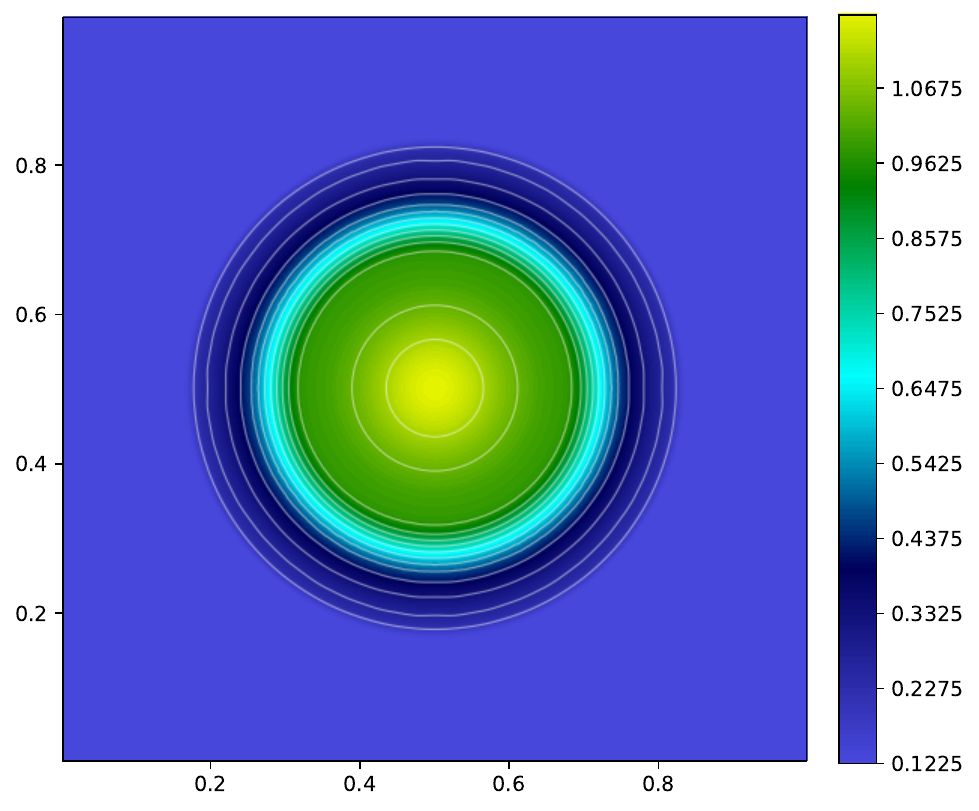}
        \caption{$t=t_f / 3$}
        \label{fig:2D_implosion_2}
    \end{subfigure}

    \bigskip

    \begin{subfigure}{0.45\textwidth}
        \includegraphics[width=\textwidth]{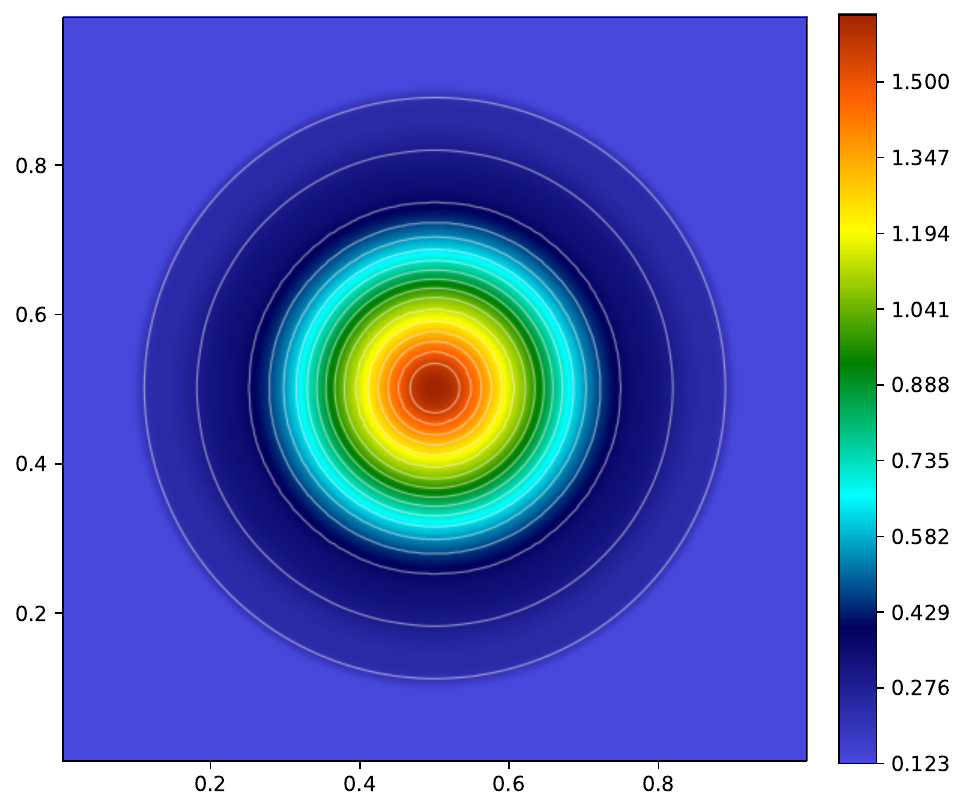}
        \caption{$t= 2 t_f / 3$}
        \label{fig:2D_implosion_3}
    \end{subfigure}
    \hfill
    \begin{subfigure}{0.45\textwidth}
        \includegraphics[width=\textwidth]{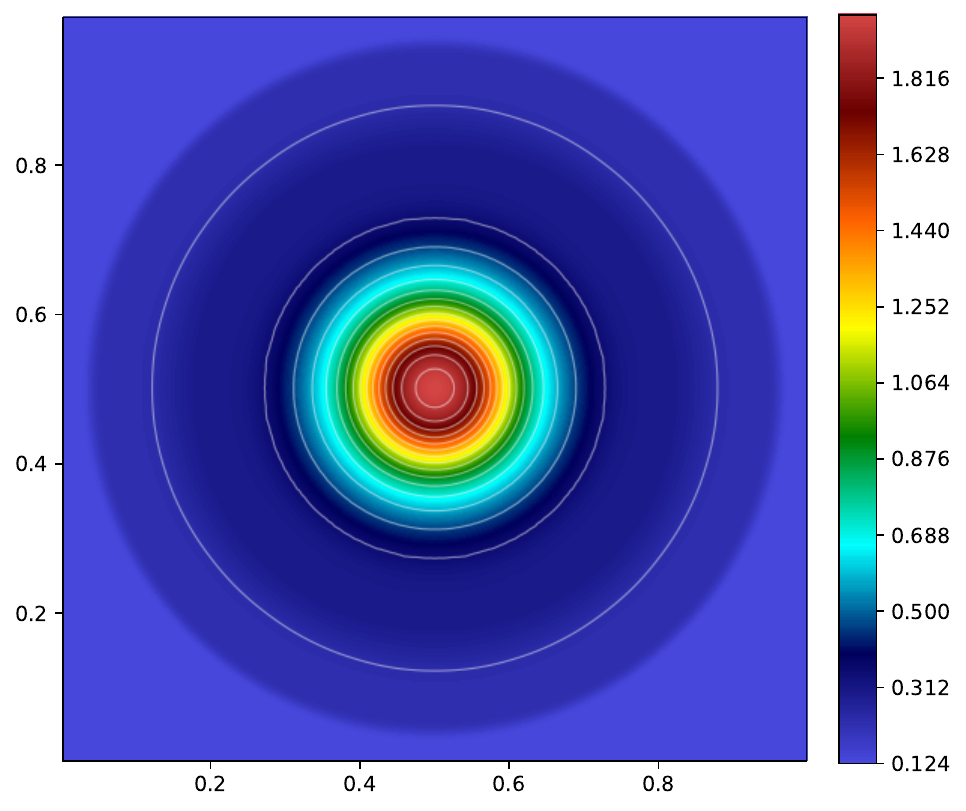}
        \caption{$t=t_f$}
        \label{fig:2D_implosion_4}
    \end{subfigure}

    \caption{Implosion test case from \cref{sec:implosion}: density $\rho$ at different times, using a mesh made of $256^2$ cells. Note that the color scale is the same for all figures.}
    \label{fig:2D_implosion}
\end{figure}

The computational domain is the unit square $[0,1]^2$,
and we take the final time $t_f = 0.125$ shortly before the acoustic waves reach the boundaries of the computational domain.
The results for the density are displayed on a mesh made of $256^2$ cells in \cref{fig:2D_implosion}, at different times.
We observe that the profile of the implosion is well-captured
by the new well-balanced scheme.

\section{Conclusions}
\label{sec:Conclusion}
In this paper, we derived a new numerical scheme to
approximate weak solutions of the Euler equations with a gravitational source term.
The presented Godunov-type scheme is based on an approximate Riemann solver
composed of two intermediate states.
It is constructed such that the scheme is consistent, positivity-preserving,
i.e., the positivity of density and pressure is preserved,
and satisfies the discrete entropy inequalities.
The proof of these properties are summarized in \cref{theo:summary} which relies on a main result of this work, given in \cref{th:entrop},
where we showed that the specific entropy
given by the HLL Riemann solver~\cite{HarLaxLee1983} satisfies a decreasing principle.
The remaining degrees of freedom in the approximate Riemann solver were exploited,
on the one hand, to incorporate the fully well-balanced property in the Riemann solver,
and on the other hand, to make the Riemann solver coincide
with the well-known HLL Riemann solver in the absence of a gravitational field.
Thus, besides the above-mentioned properties,
the new numerical scheme provably is fully well-balanced,
i.e., the scheme preserves moving equilibria and
the associated isentropic hydrostatic equilibria up to machine precision,
as proven in \cref{theo:summary}.
Further, we have shown a strategy allowing an extension of the
first-order Godunov-type numerical scheme to higher order accuracy,
while preserving the fully well-balanced property of the first-order scheme.
The performed numerical test cases illustrated
the theoretical results presented in this work.
Finally, an extension to two space dimensions is presented, where in a Cartesian framework, the well-balancing procedure is applied direction by direction. This allowed capturing two-dimensional moving equilibria, where the gravity acts only in one spatial direction.
An extension to truly two-dimensional moving equilibria has to be derived to allow the well-balancing of general multidimensional equilibria. This is closely connected to the divergence-free property, and is subject to future work.

\section*{Acknowledgments}

V. Michel-Dansac extends his thanks to ANR-22-CE25-0017 OptiTrust.
The SHARK-FV conference has greatly contributed to this work.

\bibliographystyle{plain}
\bibliography{literature}

\end{document}